\documentclass{amsart}

\makeatletter
\def\re@DeclareMathSymbol#1#2#3#4{%
    \let#1=\undefined
    \DeclareMathSymbol{#1}{#2}{#3}{#4}}
\DeclareSymbolFont{symbolsC}{U}{txsyc}{m}{n}
\SetSymbolFont{symbolsC}{bold}{U}{txsyc}{bx}{n}
\DeclareFontSubstitution{U}{txsyc}{m}{n}
\re@DeclareMathSymbol{\nsimeq}{\mathrel}{symbolsC}{59}
\makeatother

\title[Cluster Categories]{Cluster Categories}

\author[Idun Reiten]{Idun Reiten}

\index{Reiten, Idun}
\addcontentsline{toc}{author}{Idun Reiten}%
\addcontentsline{toc}{chapter}{Cluster Categories}%
\usepackage{mathrsfs}
\usepackage{bbm}
\usepackage{tikz}
\usepackage{amsmath}
\usepackage{amssymb}
\tikzstyle{dot} = [inner sep=0pt,thick,fill=black,circle,minimum size=2.5pt]
\tikzstyle{line} = [draw, -latex]

\newtheorem{theorem}{Theorem}[section]
\newtheorem{corollary}[theorem]{Corollary}

\newtheorem{example}[theorem]{Example}

\renewcommand{\a}{\alpha}
\renewcommand{\b}{\beta}
\newcommand{\g}{\gamma}
\newcommand{\fa}{f_\alpha}
\newcommand{\fb}{f_\beta}
\newcommand{\s}{\small }
\newcommand{\e}{\mathbf{e}}
\newcommand{\m}{\mbox{mod}}
\newcommand{\mkq}{\m kQ}
\newcommand{\Db}{\mathcal{D}^b(Q)}
\newcommand{\Dbk}{\mathcal{D}^b(kQ)}
\newcommand{\rep}{\mbox{rep}}
\newcommand{\CY}{$2$-CY }
\newcommand{\C}{\mathscr{C}}
\newcommand{\CQ}{\mathscr{C}_Q}
\newcommand{\Sub}{\sSub\Lambda_w}
\newcommand{\Subu}{\underline{\uSub}\Lambda_w}
\newcommand{\D}[2][b]{\mathcal{D}^{#1}(#2)}
\newcommand{\CC}[1]{\mathscr{C}_{#1}}
\newcommand{\Endop}[2][]{\End_{#1}(#2)^{\op}}%%% op changed to \op
\newcommand{\EndopC}[1]{\End_{\C}(#1)^{\op}}%%% op changed to \op
\newcommand{\EndopCQ}[1]{\End_{\CQ}(#1)^{\op}}%%% op changed to \op
\newcommand{\seq}[2][n]{\lbrace #2_{1},\ldots,\,#2_{#1} \rbrace}
\newcommand{\seqn}[2][n]{(#2_{1},\ldots,\,#2_{#1})}
\newcommand{\seqs}[2][n]{#2_{1},\ldots,\,#2_{#1}}
\newcommand{\dsum}[2][n]{#2_{1}\oplus\ldots\oplus #2_{#1}}
\newcommand{\num}[2][1]{#1,\ldots,#2}
\newcommand{\map}[2]{#1\rightarrow #2}
\newcommand{\ts}[1]{\stackrel{#1}{\rightarrow}}
\providecommand{\bysame}{\leavevmode\hbox to3em{\hrulefill}\thinspace}
\providecommand{\MR}{\relax\ifhmode\unskip\space\fi MR }
\providecommand{\href}[2]{#2}

\newcommand{\CM}{\operatorname{CM}\nolimits}
\renewcommand{\rep}{\operatorname{rep}\nolimits}
\newcommand{\Ker}{\operatorname{Ker}\nolimits}
\newcommand{\Jac}{\operatorname{Jac}\nolimits}

\newcommand{\ext}{\operatorname{Ext}\nolimits}

\renewcommand{\m}{\operatorname{mod}\nolimits}

\newcommand{\Hom}{\operatorname{Hom}\nolimits}
\newcommand{\End}{\operatorname{End}\nolimits}
\newcommand{\add}{\operatorname{add}\nolimits}

\newcommand{\sSub}{\operatorname{Sub}\nolimits}
\newcommand{\uSub}{\operatorname{Sub}\nolimits}

\newcommand{\op}{\operatorname{op}\nolimits}

\begin{document}

\begin{abstract}
Cluster algebras were introduced by Fomin-Zelevinsky in 2002 in order to give
a combinatorial framework for phenomena occurring in the context of algebraic groups. 
Cluster algebras also have links to a wide range of other subjects, including the
representation theory of finite dimensional algebras, as first discovered by Marsh-
Reineke-Zelevinsky. Modifying module categories over hereditary algebras, 
cluster categories were introduced in work with Buan-Marsh-Reineke-Todorov in
order to ``categorify'' the essential ingredients in the definition of cluster algebras
in the acyclic case. They were shown to be triangulated by Keller. Related work
was done by Geiss-Leclerc-Schr\"oer using preprojective algebras of Dynkin type. In
work by many authors there have been further developments, leading to feedback
to cluster algebras, new interesting classes of finite dimensional algebras, and the
investigation of categories of Calabi-Yau dimension $2.$

\end{abstract}

%\begin{classification}
%Primary 16G20, 13F60, 18E30; Secondary 16D90, 16G70, 05E10.
%\end{classification}

% Any keywords %
%\begin{keywords}
%quivers, path algebras, cluster algebras, cluster categories, Calabi-Yau categories, derived categories, tilting modules, almost split sequences
%\end{keywords}

% Do not remove next line %
\maketitle
\section*{Introduction}
Almost 10 years ago Fomin and Zelevinsky introduced the concept of cluster algebras, in order to create a combinatorial framework for the study of canonical bases in quantum groups, and for the study of total positivity for algebraic groups. In a series of papers they  developed a theory  of cluster algebras, which has turned out to have numerous applications to many areas of mathematics. One of the most important and influential connections has been with the representation theory of finite dimensional algebras. Such a connection was suggested by the paper \cite{MRZ}.

A cluster algebra, in its simplest form, is defined as follows, as a subalgebra of the rational function field $F =\mathbbm{Q}(x_1,\ldots, x_n).$ Start with a \emph{seed} $(\seqn{x}, B),$ which, by definition, is a pair consisting of a free generating set of $F,$ which for simplicity we choose to be $\seqn{x},$ and a skew symmetric $n\times n$ matrix $B$ over $\mathbbm{Z}.$ Alternatively, we can instead of the matrix use a finite quiver $Q$ (that is, directed graph) with  vertices $1,2,\ldots,n,$ and no oriented cycles of length $1$ or $2.$ For each $i = 1,\ldots,n,$ a new seed $\mu_i(\seqn{x},Q)$ is defined by first replacing $x_i$ with another element $x_i^*$ in $F$ according to a specific rule which depends upon both $\seqn{x}$ and $Q.$ Then we get a new free generating set $(x_1,\ldots,x_{i-1}, x_i^*, x_{i+1},\ldots, x_n).$ For $i = 1,\ldots,n$ there is defined a \textit{mutation} $\mu_i(Q)$ of the quiver $Q,$ giving a new quiver with $n$ vertices.  Then we get a new seed $((x_1,\ldots,x_{i-1}, x_i^*, x_{i+1},\ldots, x_n),\mu_i(Q)).$ We continue applying $\mu_1,\ldots, \mu_n$ to the new seeds to get further seeds. The $n$-element subsets occurring in seeds are called \textit{clusters}, and the elements in the clusters are called \textit{cluster variables.} The associated \textit{cluster algebra} is the subalgebra of $F$ generated by all cluster variables.

There are many challenging problems concerning cluster algebras. One way of attacking them is via \textit{categorification.} This is not a well defined procedure, but expresses the philosophy that we want to replace ingredients in the definition of cluster algebras  by similar concepts in a category with additional structure. The category could for example be the category of finite dimensional modules over a finite dimensional $k$-algebra for a field $k,$ or a closely related category with similar properties. In particular, the category should have enough structure so that each object $X$ has an associated finite quiver $Q_X$ (namely the quiver of the endomorphism algebra of the object). It is an extra bonus if there is a way of constructing the original ingredients of the cluster algebra back from the chosen analogous object in the new category, but we do not require that this should always be the case.

We first discuss the special case of categorifying quiver mutation alone. Then we want to find a ``nice'' category $\C$ with a distinguished set $\mathscr{T}$ of objects, with an operation $T\mapsto\mu_i(T)$ defined for $T$ in $\mathscr{T}$ and any $i= 1,\ldots, n.$ We would like this operation to ``lift'' the quiver mutation, that is, we would want that $Q_T = Q$ and $Q_{\mu_i(T)} = \mu_i(Q).$ 

We then discuss categorification of some of the essential ingredients involved in the definition of cluster algebras, such as clusters, cluster variables, seeds. We want to imitate these concepts, and preferably also operations of addition and multiplication involving them, in a ``nice'' category $\C.$ As analogs of clusters we want a distinguished set $\mathscr{T}$ of objects of the form $T = \dsum{T},$ where the $T_i$ are indecomposable and $T_i\nsimeq T_j$ for $i\neq j.$ The $T_i$ would then be the analogs of cluster variables. For each $i= 1,\ldots,n$ we want a unique indecomposable object $T_i^*\nsimeq T_i,$ where $T_i^*$ is a summand of an object in $\mathscr{T},$ such that $T/T_i\oplus T_i^*$ is in $\mathscr{T}.$ To find analogs of the seeds, we consider pairs $(T,Q)$ where $T$ is in $\mathscr{T}$ and $Q$ is a quiver with $n$ vertices. For a seed $(\seqn{u},Q)$ it does not make sense to talk about a connection between $\seqn{u}$ and $Q.$ But for the pair $(T,Q)$ a natural connection to ask for between $T$ and $Q$ is that $Q$ coincides with the quiver $Q_T$ of the endomorphism algebra $\Endop{T}.$ We can try to choose $T$ such that we have an \textit{initial tilting seed} $(T,Q)$ with this property. If the same set $\mathscr{T}$ provides a categorification of quiver mutation in the sense discussed above, then this nice property for a tilting seed will hold for all pairs obtained from $(T,Q)$ via a sequence of mutations of the objects\break in $\mathscr{T}.$ 

%Here we explain how we can find an appropriate categorification in
%the case of acyclic cluster algebras. We let $Q$ be an acyclic quiver with $n$ vertices.

We have collected a (not complete) list of desired properties for the
categories, together with a distinguished set of objects $\mathscr{T},$ which we would like our categorification to satisfy. But there is of course no guarantee to start with that it is possible to find a satisfactory solution. Here we explain how we can find an appropriate
categorification in the case of what is called acyclic cluster
algebras. We let $Q$ be an acyclic quiver, that is, a quiver
with no oriented cycles, with $n$ vertices. Given a field $k,$ there is a way of associating a finite dimensional path algebra $kQ$ with the quiver $Q$ (see Section \ref{sec1}). Then the category of finite dimensional $kQ$-modules might be a candidate for the category we are looking for. The tilting modules, which have played a central role in the representation theory of finite dimensional algebras, might be a candidate for the distinguished set of objects, since they have some of the desired properties, as also suggested by the work in \cite{MRZ}. A $kQ$-module $T=\dsum{T}$ is a \textit{tilting} module if the $T_i$ are indecomposable, $T_i\nsimeq T_j$ for $i\neq j,$ and every exact sequence of the form $0\rightarrow T\rightarrow E\rightarrow T\rightarrow 0$ splits. However, this choice does not quite work since it may happen that for some $i$ there is no indecomposable module $T_i^*\nsimeq T_i$ such that $T/T_i\oplus T_i^*$ is a tilting module. The idea is then to ``enlarge'' the category $\mkq$ to make it more likely to find some $T_i^*.$
 This enlargement can in practice be done by taking a much larger category containing $\mkq,$ namely the bounded derived category $\Dbk,$ and then taking the orbit category under the action of a suitable cyclic group in order to cut down the size. Then we end up with what has been called the \textit{cluster category} $\CC{Q}$ \cite{BMRRT}. As distinguished set of objects $\mathscr{T}$ we choose an enlargement of the set of tilting $kQ$-modules, called \textit{cluster tilting} objects. Then $\CC{Q},$ together with $\mathscr{T},$ has all the properties we asked for above, and some more which we did not list.
 
It is natural to try to find other categories with distinguished sets of objects which would categorify other classes of cluster algebras. This has been succesfully done in \cite{GLS1} using the category of finite dimensional modules $\m\Lambda$ over a class of finite dimensional algebras $\Lambda$ called preprojective algebras of Dynkin type. 

The investigation of cluster categories and preprojective algebras of Dynkin type has further led to work on what is called Hom-finite triangulated $2$-Calabi-Yau categories (\CY for short), with a specific set of objects called cluster tilting objects \cite{IY}. The endomorphism algebras of the cluster tilting objects form an interesting class of finite dimensional algebras. For example, the investigation of the cluster tilted algebras, which by definition are those coming from cluster categories, has shed new light on tilting theory.

In addition, there have been exciting applications of the categorification to the theory of cluster algebras. For this it is useful to establish a tighter \hbox{connection} between cluster algebras and cluster categories (or more general \CY categories). This can happen through providing maps in one or both directions, especially between the cluster variables and the indecomposable summands of the cluster tilting objects. In many cases, for example for cluster categories, such explicit maps have been constructed, giving deep connections which are useful for applications to cluster algebras (\cite{BMRT},\cite{CC},\cite{CK1},\cite{CK2},\cite{GLS1},\cite{P}).

A special case of the problem of categorifying quiver mutation was investigated in the early days of the present form of the representation theory of finite dimensional algebras, which started around 40 years ago. Reflections of quivers at vertices which are sinks (or sources) were introduced in \cite{BGP}. A \textit{sink} $i$ is a vertex where no arrow starts, and a \textit{source} is a vertex where no arrow ends. The operation $\mu_i$ on a quiver where $i$ is a sink was defined by reversing all arrows ending at $i.$ A categorification of this special case of quiver mutation was done using tilting modules over the path algebras $kQ$ for a finite acyclic quiver $Q$ \cite{APR}. In this special case it was possible to use tilting\break modules \cite{APR}. 

We start the paper with a discussion of this mutation from \cite{BGP} and its categorification in Section \ref{sec1}. We also give some definitions and basic properties of path algebras and quiver representations. In Section \ref{sec2} we give an introduction to the theory of cluster algebras, including definitions, examples and crucial properties. Cluster categories are introduced in Section \ref{sec3}. We give some motivation, including a list of desired properties which they should satisfy, and illustrate through examples. We also introduce the cluster tilting objects. In Section \ref{sec4} we deal with generalizations to Hom-finite  triangulated \CY categories. The endomorphism algebras of cluster tilting objects are called \textit{\CY-tilted algebras.} They are discussed in Section \ref{sec5}, together with their relationship to the interesting class of Jacobian algebras, which are given by quivers with potential \cite{DWZ1}. In Section \ref{sec6} we discuss applications to cluster algebras in the acyclic case. 

The various aspects of the relationship between cluster algebras and the representation theory of finite dimensional algebras have stimulated a lot of research activity during the last few years, with several interesting developments by a large number of contributors. I have chosen to emphasize aspects closest to my own interest, which deal with the more categorical aspects of the subject. But cluster algebras lie in the center of it all, as inspiration, so I have included a brief discussion of them, as well as a discussion of the  feedback of the general categorical approach to cluster algebras.

Several important topics related to cluster algebras are not discussed in this paper. In particular, this concerns the developments dealing with the cluster algebras themselves, the applications of categorification to the construction of (semi)canonical bases and their duals \cite{GLS1}, and the recent work in \cite{KS} related to $3$-Calabi-Yau algebras. Another interesting type of categorification was \hbox{introduced} in \cite{HL}, and pursued in \cite{N}, see also \cite{L}. We also refer to \cite{DWZ2} for work using quivers with potentials to solve a series of conjectures.

We refer to the  survey papers (\cite{BM},\cite{Fo},\cite{K3},\cite{K5},\cite{K6},\cite{K7},\cite{L},\cite{Re1},\cite{Re2},\cite{Ri}, \linebreak \cite{Ze}) for additional information.
 
\section{Bernstein-Gelfand-Ponomarev Reflections}\label{sec1}
In this section we go back to the beginning of the present form of the representation theory of finite dimensional algebras, which  dates back to around 1970. We show that some of the early developments can be seen as categorification of Bernstein-Gelfand-Ponomarev reflections of quivers at sinks (or sources), which is a special case of the quiver mutation used by Fomin-Zelevinsky in connection with their definition of cluster algebras. This categorification involves tilting theory (\cite{BB},\cite{HR}), which is one of the most important developments in the representation theory of algebras. There are numerous applications both within the field and outside. We start with some relevant background material. We refer to the books (\cite{ASS},\cite{ARS},\cite{GR},\cite{Ha2},\cite{Ri2}) as general references.

%Let $Q$ be a finite quiver, that is, a quiver with a finite number of
%vertices and a finite number of arrows between them. 

\subsection{Representations of quivers}
Let $Q$ be a {\it finite quiver}, that is, a directed graph
with a finite number of vertices and a finite number of arrows
between the vertices.
Assume that the quiver $Q$ is {\it acyclic}, that is, has
no oriented cycles; for example, let $Q$ be the quiver\begin{minipage}{2.3cm}
\begin{center} 
\begin{tikzpicture}
\node[] (A) at (0,0) {\s $1$};
\node[] (B) at (1,0) {\s $2$};
\node[] (C) at (2,0) {\s $3$};
\node[] () at (0.5,0.2) {\s $\a$};
\node[] () at (1.5,0.2) {\s $\b$};
\draw[->] (0.1,0) -- (0.9,0);
\draw[->] (1.1,0) -- (1.9,0);
\end{tikzpicture}
\end{center}
\end{minipage}. For simplicity of exposition, we define the most relevant concepts only for this example.

Let $k$ be a field, which we always assume to be algebraically closed.
A representation of $Q$ over $k$  is $V_1\stackrel{\fa}{\rightarrow}V_2\stackrel{\fb}{\rightarrow}V_3,$ where we have associated a finite dimensional vector space $V_i$ to each vertex $i,$ and a linear transformation to each arrow. A map
$$h:(V_1\stackrel{\fa}{\rightarrow}V_2\stackrel{\fb}{\rightarrow}V_3)\longrightarrow(V_1'\stackrel{\fa'}{\rightarrow}V_2'\stackrel{\fb'}{\rightarrow}V_3')$$
between two representations is a triple $h = (h_1,h_2,h_3),$ where $h_i:V_i\rightarrow V_i'$ is a linear transformation for each $i,$ such that the following diagram commutes.
\begin{center} 
\begin{tikzpicture}
\node[] () at (0,0) {$V_1$};
\node[] () at (2,0) {$V_2$};
\node[] () at (4,0) {$V_3$};
\node[] () at (0,-1) {$V_1'$};
\node[] () at (2,-1) {$V_2'$};
\node[] () at (4,-1) {$V_3'$};
\draw[->] (0.2,0) -- (1.8,0);
\draw[->] (2.2,0) -- (3.8,0);
\draw[->] (0.2,-1) -- (1.8,-1);
\draw[->] (2.2,-1) -- (3.8,-1);
\draw[->] (0,-0.2) -- (0,-0.8);
\draw[->] (2,-0.2) -- (2,-0.8);
\draw[->] (4,-0.2) -- (4,-0.8);
\node[] () at (0.2,-0.5) {\s $h_1$};
\node[] () at (2.2,-0.5) {\s $h_2$};
\node[] () at (4.2,-0.5) {\s $h_3$};
\node[] () at (1,0.2) {\s $\fa$};
\node[] () at (3,0.2) {\s $\fb$};
\node[] () at (1,-0.8) {\s $\fa'$};
\node[] () at (3,-0.8) {\s $\fb'$};
\end{tikzpicture}
\end{center}

The category $\rep Q$ of representations of $Q,$ with objects and maps as defined above, is equivalent to the category $\m kQ$ of finite dimensional modules over the path algebra $kQ.$ Here the paths in $Q,$ including the trivial paths $e_i$ associated with the vertices $i,$ are a $k$-basis for $kQ.$  So in our example $\{\a,\, \b,\, \b\a,\, e_1,\, e_2,\, e_3\}$ is a $k$-basis.
The multiplication for the basis elements is defined as composition of paths whenever this is possible, and is defined to be $0$ otherwise. For example, we have \sloppy $\b\cdot\a = \b\a,\,\, \a\cdot e_2 = 0,\,\, \a\cdot e_1 = \a,\break\,e_1\cdot e_1 = e_1.$

The connection between $\rep Q$ and $\m kQ$ is illustrated as
follows.  If
$V_1\stackrel{\fa}{\rightarrow}V_2\stackrel{\fb}{\rightarrow}V_3$ is in $\rep Q,$ then the vector space $V_1\oplus V_2\oplus V_3$ 
can be given a $kQ$-module structure by defining $\a(v_1,v_2,v_3) = (0,\fa(v_1),0),\,\, \b(v_1,v_2,v_3) =
 (0,0,\fb(v_2)),\, e_1(v_1,v_2,v_3) \kern-.2pt= \kern-.2pt(v_1,0,0),\, e_2(v_1,v_2,v_3) \kern-.2pt= \kern-.2pt(0,v_2,0),\, 
e_3(v_1,v_2,v_3) \kern-.2pt= (0,0,v_3).$

The indecomposable projective representations $P_1,\,P_2$ and $P_3$ associated with the vertices $1,2$ and $3$ are $k\stackrel{id}{\rightarrow}k\stackrel{id}{\rightarrow}k,\,\,0\rightarrow k\stackrel{id}{\rightarrow}k$ and $0\rightarrow 0\rightarrow k,$ the simple representations $S_1,\,S_2$ and $S_3$ are $k\rightarrow 0\rightarrow 0,\,\,0\rightarrow k\rightarrow 0$ and $0\rightarrow 0\rightarrow k$ and the indecomposable injective representations $I_1, I_2$ and $I_3$ are $k\rightarrow 0\rightarrow 0,\, k\stackrel{id}{\rightarrow}k\rightarrow 0$ and $k\stackrel{id}{\rightarrow}k\stackrel{id}{\rightarrow}k.$ We also use the same notation for a representation viewed as a $kQ$-module.

An important early result on quiver representations was the following \cite{Ga}.

\begin{theorem}{\normalfont}\label{thm11}
Let $Q$ be a finite connected quiver and $k$ an algebraically closed field. Then $\rep Q$ has only a finite number of indecomposable representations up to isomorphism if and only if the underlying graph $|Q|$ is of Dynkin type $A_n, D_n, E_6, E_7$ or $E_8.$
\end{theorem}

\subsection{Reflection functors}\label{refunc}
The proof by Gabriel of Theorem \ref{thm11} was technically complicated. A more elegant proof was soon thereafter given by Bernstein-Gelfand-Ponomarev \cite{BGP}, taking advantage of the fact that the classification theorem involved Dynkin diagrams, a fact which suggested connections with root systems and positive definite quadratic forms. One important aspect of their work was that they introduced reflections of quivers and associated reflection functors. Using this, together with some special modules from \cite{APR} inspired by the reflection functors, later known as APR-tilting modules, we shall illustrate the idea of categorification.  These examples are special cases of categorifications of more general mutation of quivers, which is important for the categorification of cluster algebras.

Let $i$ be a vertex in the quiver $Q$ which is a sink, that is, there are no arrows starting at $i.$ In our running example the vertex $3$ is a sink. We define a new quiver $\mu_3(Q),$ known as the \textit{mutation} of $Q$ at the vertex $3.$ This is obtained by reversing all the arrows in $Q$ ending at $3,$ so in our example it is the\break quiver\begin{minipage}{2.5cm}
\begin{center} 
\begin{tikzpicture}
\node[] (A) at (0,0) {\s $1$};
\node[] (B) at (1,0) {\s $2$};
\node[] (C) at (2,0) {\s $3.$};
\draw[->] (0.1,0) -- (0.9,0);
\draw[<-] (1.1,0) -- (1.9,0);
\end{tikzpicture}
\end{center}
\end{minipage} In \cite{BGP} a reflection functor $F_3:\rep Q\rightarrow\rep Q'$ was defined on objects by sending $V_1\stackrel{\fa}{\rightarrow}V_2\stackrel{\fb}{\rightarrow}V_3$ to $V_1\stackrel{\fa}{\rightarrow}V_2\stackrel{\fb'}{\leftarrow}\Ker\fb,$ where $\fb'$ is the  natural inclusion. We see that $S_3$ is sent to the zero representation, and we have the following connection between $\rep Q$ and $\rep Q',$ where $S_3'$ is the representation $0\rightarrow 0\leftarrow k$ of $Q'.$

\begin{theorem}\label{Thm12}
The reflection functor $F_3$ induces an equivalence $F_3:\rep Q\backslash S_3\rightarrow\rep Q'\backslash S_3',$ where $\rep Q\backslash S_3$ denotes the full subcategory of $\rep Q$ \hbox{consisting} of  objects which are finite direct sums of indecomposable objects not isomorphic to $S_3.$
\end{theorem}

This is a key step in the proof in \cite{BGP} of Gabriel's theorem. It gives an easy illustration of a categorification, where a quiver $Q$ is replaced  by some object in $\rep Q,$ and the mutation of $Q$ at a sink $i$ is replaced by the associated functor~$F_i.$

\subsection{Illustration using AR-quiver}\label{Arquiv}
For a finite dimensional $k$-algebra $\Lambda$ a special kind of exact sequence, known as an \textit{almost split sequence} (or also \textit{Auslander-Reiten} sequence), was introduced in \cite{AR}. An exact sequence $0\rightarrow A\stackrel{f}{\rightarrow}B\stackrel{g}{\rightarrow}C\rightarrow 0$ is almost split if it is not split, the end terms are indecomposable, and each map $h:\map{X}{C}$ where $X$ is indecomposable and $h$ is not an isomorphism, factors through $g:\map{B}{C}.$  We have the following basic result, where we assume that all our modules are finite dimensional over\break $k$ \cite{AR}.

\begin{theorem}
For any indecomposable nonprojective $\Lambda$-module $C$ (or for any indecomposable noninjective $\Lambda$-module $A$), there exists an almost split sequence $0\rightarrow A\rightarrow B\rightarrow C\rightarrow 0,$ which is unique up to isomorphism.
\end{theorem}

The almost split sequences induce an operation $\tau,$ called the AR-translation, from the indecomposable nonprojective $\Lambda$-modules to the indecomposable noninjective ones, satisfying $\tau(C) = A$ when $0\rightarrow A\rightarrow B\rightarrow C\rightarrow 0$ is almost split.

On the basis of the information given by the almost split sequences (in general together with some special maps to projectives and from injectives) we can draw a new quiver, called the \textit{Auslander-Reiten quiver} (AR-quiver for short), where the vertices correspond to the isomorphism classes of indecomposable $\Lambda$-modules. 

In our examples we have the following.
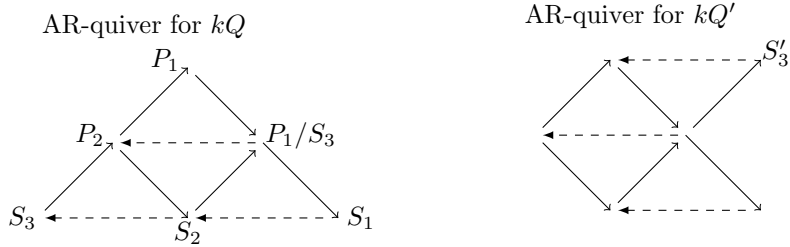
\newline~\newline
\begin{minipage}{6.5cm}
\begin{center}
AR-quiver for $kQ$\newline
\begin{tikzpicture}
\draw[->] (0.1,0.1) -- (1,1);
\draw[->] (1.1,1.1) -- (2,2);
\draw[->] (1.1,0.9) -- (2,0);
\draw[->] (2.1,1.9) -- (2.9,1.1);
\draw[->] (2.1,0.1) -- (2.9,0.9);
\draw[->] (3,1) -- (4,0);
\path[line,dashed] (1.8,0) -- (0.1,0);
\path[line,dashed] (3.8,0) -- (2.1,0);
\path[line,dashed] (2.9,1) -- (1.1,1);
\node[] () at (-0.2,0) {$S_3$};
\node[] () at (0.7,1.1) {$P_2$};
\node[] () at (1.7,2.1) {$P_1$};
\node[] () at (2,-0.2) {$S_2$};
\node[] () at (4.3,0) {$S_1$};
\node[] () at (3.5,1.1) {$P_1/S_3$};
\end{tikzpicture}
\end{center}
\end{minipage}
\begin{minipage}{6cm}
\begin{center}
AR-quiver for $kQ'$\newline
\begin{tikzpicture}
\draw[->] (3.1,1.1) -- (4,2);
\draw[->] (1.1,1.1) -- (2,2);
\draw[->] (1.1,0.9) -- (2,0);
\draw[->] (2.1,1.9) -- (2.9,1.1);
\draw[->] (2.1,0.1) -- (2.9,0.9);
\draw[->] (3,1) -- (4,0);
\path[line,dashed] (3.8,2) -- (2.1,2);
\path[line,dashed] (3.8,0) -- (2.1,0);
\path[line,dashed] (2.9,1) -- (1.1,1);
\node[] () at (4.2,2.1) {$S_3'$};
\node[] () at (2,-0.6) {};
\end{tikzpicture}
\end{center}
\end{minipage}\newline\newline
The broken arrows indicate the translation $\tau,$ and we can then deduce the shape of the almost split sequences from the AR-quiver. For example, for $kQ$ we have the almost split sequences $0\rightarrow S_3\rightarrow P_2\rightarrow S_2\rightarrow 0,\,\, 0\rightarrow S_2\rightarrow P_1/S_3\rightarrow S_1\rightarrow 0$ and $0\rightarrow P_2\rightarrow P_1\oplus S_2\rightarrow P_1/S_3\rightarrow 0.$ The AR-quivers for $kQ$ and $kQ'$ are not isomorphic, but when dropping $S_3$ from the first one and $S_3'$ from the second one, they are clearly isomorphic. This reflects the fact that there is an equivalence between the corresponding subcategories, as stated in Theorem \ref{Thm12}.

\subsection{Module theoretical interpretation}\label{Modth}
Let $Q$ be a finite quiver without oriented cycles and with vertices $1,\ldots,n.$ Let $i$ be a vertex of $Q$ which is a sink. Denote by $\mu_i(Q) = Q'$ the quiver obtained by  mutation at $i.$ We write $kQ = P_1\oplus\ldots\oplus P_n,$ where $P_j$ is the indecomposable projective $kQ$-module associated with the vertex $j.$ Then we have the following module theoretical  interpretation of the reflection functors \cite{APR}.

\begin{theorem}
With the above notation, we have the following.
\begin{enumerate}\leftskip-2pt
	\item For $T = kQ/P_i\oplus\tau^{-1}P_i$ we have $\Endop[kQ]{T} \simeq kQ'.$
	\item The functor $F_i:\rep Q\rightarrow\rep Q'$ is isomorphic to the functor $\Hom_{kQ}(T,\,)\!:\m kQ\rightarrow\m kQ'.$
\end{enumerate}
\end{theorem}

Let $J=J_Q$ be the ideal generated by all arrows in a path algebra $kQ.$
It is known that any finite dimensional $k$-algebra $A$ is Morita equivalent to $kQ/I$ for some finite quiver $Q$ and ideal $I$ in $kQ$ with $I\subseteq J^2.$ A generating set of $I$ is called a set of \textit{relations} for $A.$  We denote $Q$ by $Q_A.$ In particular, we have $Q_{kQ} = Q$ and $Q_{kQ'} = Q'.$ If $X$ is in $\m A,$ then the associated quiver $Q_X$ is by definition the quiver $Q_{\mbox{\scriptsize End}(X)^{\mbox{\scriptsize  op}}}$ associated with the finite dimensional $k$-algebra $\End(X)^{\op}.$ 

In our running example we have  $T = P_1\oplus P_2\oplus\tau^{-1}S_3.$ Consider the diagram
\begin{center}
\begin{tikzpicture}
\draw[->] (0.2,0) -- (1.4,0);
\draw[-] (0.2,-0.1) -- (0.2,0.1);
\draw[->] (0.3,1) -- (2.5,1);
\draw[-] (0.3,0.9) -- (0.3,1.1);
\draw[->] (0,0.8) -- (0,0.2);
\draw[-] (-0.1,0.8) -- (0.1,0.8);
\draw[->] (2.7,0.8) -- (2.7,0.2);
\draw[-] (2.6,0.8) -- (2.8,0.8);
\node[] () at (0,0) {$Q$};
\node[] () at (0,1) {$kQ$};
\node[] () at (2.7,1) {$T$};
\node[] () at (2.3,0) {$\mu_3(Q)=Q'$};
\node[] () at (0.8,0.14) {\s $\mu_3$};
\end{tikzpicture}
\end{center}
Then $Q_T = Q'$ since $\Endop[kQ]{T}\simeq kQ'.$ Further, we can define $\mu_3(kQ)$ by replacing $P_3 = S_3$ by $\tau^{-1}S_3,$ so that $\mu_3(kQ)=T.$  This way we can view the above theorem as a way of categorifying quiver mutation. This categorification, which is quite different from the one discussed in Section \ref{refunc}, is of the type we shall be dealing with, and it will be generalized later.

\subsection{Tilting theory}\label{tilttheory}
The $kQ$-modules $kQ$ and $T$ are examples of what are now called tilting modules  (\cite{BB}, \cite{HR}). A module $T$ over a path algebra $kQ$ is a \emph{tilting module} if  $\ext_{kQ}^1(T,T) = 0,$ and the number of nonisomorphic indecomposable summands is the number of vertices in the quiver $Q$. Objects $X$ with $\ext_{kQ}^1(X,X)= 0$ are called \textit{rigid.} The endomorphism algebras $\Endop[kQ]{T}$ are by definition the \emph{tilted} algebras.
It was the work discussed in Sections \ref{refunc} and \ref{Modth} which inspired tilting theory.

An aspect of  tilting theory of interest in this paper is the following (\cite{RS},\cite{U},\cite{HU1}).

\begin{theorem}
Let $T=T_1\oplus\ldots \oplus T_n$ be a tilting $kQ$-module, where the $T_i$ are indecomposable and $T_i\nsimeq T_j$ for $i\neq j.$
\begin{enumerate}
\item For each $i,$ there is at most one indecomposable $kQ$-module $T_i^*\nsimeq T_i$ such that $T/T_i\oplus T_i^*$ is a tilting $kQ$-module.
\item For each $i$ there exists such a module $T_i^*$ if and only if $T/T_i$ is a sincere $kQ$-module, that is, all simple $kQ$-modules are composition factors of $T/T_i.$
\end{enumerate}
\end{theorem}

The modules $T/T_i$ are called \textit{almost complete} tilting modules, and $T_i$ (and $T_i^*$ if it exists) are called \textit{complements} of $T/T_i.$
The $kQ$-module $kQ$ is clearly a tilting module, and when $i$ is a sink in the quiver $Q,$ then there is always some indecomposable $kQ$-module $P_i^*\nsimeq P_i$ such that $kQ/P_i\oplus P_i^*$ is a tilting module. The module $T_i^*$ is in fact $\tau^{-1}P_i.$ So we can view the previously defined operation $\mu_i(kQ)$ as replacing $P_i$ by the unique indecomposable $kQ$-module $P_i^*\nsimeq P_i$ such that $kQ/P_i\oplus P_i^*$ is a tilting module.

\subsection{General quiver mutation and tilting modules}\label{genquivmut}
For quivers with no loops 
\begin{minipage}{0.5cm}
\begin{center} 
\begin{tikzpicture}
\node[dot] () at (0,0) {};
\draw[->] (0,0.1) arc (0:310:0.2);
\end{tikzpicture}
\end{center}
\end{minipage}
and no (oriented) $2$-cycles 
\begin{minipage}{1.2cm}
\begin{center} 
\begin{tikzpicture}
\node[dot] () at (0,0.2) {};
\node[dot] () at (0.6,0.2) {};
\draw[->] (0.1,0.2) -- (0.5,0.2);
\draw[<-] (0,0) -- (0.6,0);
\end{tikzpicture}
\end{center}
\end{minipage},
Fomin and Zelevinsky have  introduced a mutation of quivers at any vertex of the quiver as follows \cite{FZ1} (see also papers in mathematical physics \cite{Sei}). 

Let $i$ be a vertex in the quiver $Q.$
\begin{itemize}
\item[(i)] Each pair of arrows $s\rightarrow i\rightarrow t$ in $Q$ gives rise to a new arrow $s\rightarrow t$ in the mutated quiver $\mu_i(Q).$
\item[(ii)] We reverse the arrows starting or ending at $i.$
\item[(iii)] We remove any $2$-cycles.
\end{itemize}

\begin{example}\normalfont
Let $Q$ be the quiver 
\begin{minipage}{2.4cm}
\begin{center} 
\begin{tikzpicture}
\node[] (A) at (0,0) {\s $1$};
\node[] (B) at (1,0) {\s $2$};
\node[] (C) at (2,0) {\s $3$};
\node[] () at (0.5,0.2) {\s $\a$};
\node[] () at (1.5,0.2) {\s $\b$};
\draw[->] (0.1,0) -- (0.9,0);
\draw[->] (1.1,0) -- (1.9,0);
\end{tikzpicture}
\end{center}
\end{minipage}. Then $\mu_2(Q)$ is the quiver \newline\hspace*{-6pt}
\begin{minipage}{2.4cm}
\begin{center} 
\begin{tikzpicture}
\node[] (A) at (0,0) {\s $1$};
\node[] (B) at (1,0) {\s $2$};
\node[] (C) at (2,0) {\s $3$};
\draw[->] (0.1,0.1) arc (120:60:1.8);
\draw[<-] (0.1,0) -- (0.9,0);
\draw[<-] (1.1,0) -- (1.9,0);
\end{tikzpicture}
\end{center}
\end{minipage}.
Let again $kQ=P_1\oplus P_2\oplus P_3$ be the tilting kQ-module associated with $Q.$ Then $T=P_1\oplus S_1\oplus P_3$ is also a tilting $kQ$-module, and we define $\mu_2(kQ) = T$ according to the previous principle. But the quiver $Q_T$ associated to $T$ can be shown to be 
\begin{minipage}{2.4cm}
\begin{center} 
\begin{tikzpicture}
\node[] (A) at (2,0) {\s $3$};
\node[] (A) at (0,0) {\s $1$};
\node[] (A) at (1,0) {\s $2$};
\draw[->] (0.1,0.1) arc (120:60:1.8);
\draw[<-] (0.1,0) -- (0.9,0);
\end{tikzpicture}
\end{center}
\end{minipage}. Hence mutation of tilting modules does not give a categorification of quiver mutation in this case.
For the vertex $1,$ it is not even possible to replace $P_1$ to get another tilting module.
\end{example}

In conclusion, we can use tilting modules to categorify quiver mutation at sinks, but not at an arbitrary vertex. We shall see in Section \ref{sec3}  how we can modify the module category $\m kQ,$ and use objects related to tilting \hbox{$kQ$-modules}, in order to make things work for quiver  mutation at any vertex of an acyclic quiver (and more generally any quiver in the mutation class of an acyclic quiver). We shall also see that in addition to being interesting in itself, this work can be used to obtain information on the mutation class of a finite acyclic quiver. Here we take advantage of the richer structure provided by the categorification. We shall also see that for another class of quivers one can actually use tilting modules  to categorify quiver mutation \cite{IR1}.

\section{Cluster Algebras}\label{sec2}
Cluster algebras were introduced by Fomin and Zelevinsky in \cite{FZ1}.  The motivation was to create a common framework for phenomena occurring in connection with total positivity and canonical bases. The theory has had considerable influence on many different areas, amongst them the theory of quiver representations.  In this section we give basic definitions and state some main results, following  (\cite{FZ1},\cite{FZ3}), in order to have appropriate background for discussing categorification of cluster algebras.

\subsection{Cluster algebras with no coefficients}
We first discuss  cluster algebras with ``no coefficients,'' which we mainly deal with in this paper. Let $F = \mathbb{Q}\seqn{x}$ be the function field in $n$ variables over the field $\mathbb{Q}$ of rational numbers, and let $B = (b_{ij})$ be an $n\times n$ - matrix. We assume for simplicity that $B$ is skew symmetric. Then $B$ corresponds to a quiver with $n$ vertices, where there are $b_{ij}$ arrows from  $i$ to $j$ if $b_{ij}>0.$ Here we deal with quivers instead of matrices. We start with an \emph{initial seed} $(\seqn{x},Q),$ consisting of a free generating set $\seqn{x}$ for $F,$ together with a finite quiver $Q$ with $n$ vertices, labelled $\num{n}.$ For each $i \in \{\num{n}\}$ we define a new seed $\mu_i(\seqn{x}, Q)$ to be a pair $((x_1,\ldots,x_{i-1},x_i^*,x_{i+1},\ldots,x_n),\mu_i(Q)).$ Here $x_i^*$ is defined by the equality $x_ix_i^* = m_1+m_2,$ where $m_1$ is the product whose terms are $x_j^s$ if there are s arrows from $j$ to $i,$ and $m_2$ is the corresponding product associated with the arrows starting at $i.$ If $i$ is a source, that is, no arrow ends at $i,$ we set $m_1 = 1,$ and if $i$ is a sink we set $m_2 = 1.$ It can be shown that $(x_1,\ldots,x_{i-1},x_i^*,x_{i+1},\ldots,x_n)$ is again a free generating set, and that $\mu_i^2(\seqn{x},Q) = (\seqn{x},Q).$ Note that the quiver $\mu_i(Q)= Q'$ only depends on the quiver $Q,$ while the new free generating set depends on both $Q$ and the old free generating set.

We illustrate with the following.

\begin{example}\normalfont
Let $F = \mathbb{Q}(x_1,x_2,x_3)$ and $Q$ be{\vspace*{-5pt}} the quiver
\begin{minipage}{2.27cm}
\begin{center} 
\begin{tikzpicture}
\node[] () at (0,0) {\s $1$};
\node[] () at (1,0) {\s $2$};
\node[] () at (2,0) {\s $3$};
\draw[->] (0.1,0) -- (0.9,0);
\draw[->] (1.1,0) -- (1.9,0);
\end{tikzpicture}
\end{center}
\end{minipage}. 
We start with the{\vspace*{-5pt}} seed 
\begin{minipage}{3.6cm}
\begin{center} 
\begin{tikzpicture}
\node[] () at (-1.1,0) {$\Bigl((x_1,\, x_2,\, x_3),\,$};
\node[] () at (0,0) {\s $1$};
\node[] () at (0.5,0) {\s $2$};
\node[] () at (1,0) {\s $3$};
\draw[->] (0.1,0) -- (0.4,0);
\draw[->] (0.6,0) -- (0.9,0);
\node[] () at (1.2,0) {$\Bigr)$};
\end{tikzpicture}
\end{center}
\end{minipage}, and apply $\mu_3.$ Then\break\hspace*{-10pt}
\begin{minipage}{4.5cm}
\begin{center} 
\begin{tikzpicture}
\node[] () at (-0.4,0) {$\mu_3\Bigl($};
\node[] () at (0,0) {\s $1$};
\node[] () at (0.5,0) {\s $2$};
\node[] () at (1,0) {\s $3$};
\draw[->] (0.1,0) -- (0.4,0);
\draw[->] (0.6,0) -- (0.9,0);
\node[] () at (1.7,0) {$\Bigr)=Q':$};
\node[] () at (2.5,0) {\s $1$};
\node[] () at (3,0) {\s $2$};
\node[] () at (3.5,0) {\s $3.$};
\draw[->] (2.6,0) -- (2.9,0);
\draw[<-] (3.1,0) -- (3.4,0);
\end{tikzpicture}
\end{center}
\end{minipage} Furthermore{\vspace*{-5pt}} $x_3^*x_3 = x_2+1,$  so that $x_3^* = \frac{x_2+1}{x_3}.$ Hence we obtain the new seed\begin{minipage}{3.85cm}
\begin{center} 
\begin{tikzpicture}
\node[] () at (-1.25,0) {$\Bigl((x_1,x_2,\frac{x_2+1}{x_3}),\,$};
\node[] () at (0.00,0) {\s $1$};
\node[] () at (0.50,0) {\s $2$};
\node[] () at (1.0,0) {\s $3$};
\draw[->] (0.1,0) -- (0.4,0);
\draw[<-] (0.6,0) -- (0.9,0);
\node[] () at (1.2,0) {$\Bigr)$};
\end{tikzpicture}
\end{center}
\end{minipage}. Similarly, we get $\mu_2((x_1,\, x_2,\, x_3),Q) =$\,\begin{minipage}{3.95cm}
\begin{center} 
\begin{tikzpicture}
\node[] () at (-1.33,0) {$\Bigl((x_1,\frac{x_1+x_3}{x_2},x_3),\,$};
\node[] () at (0,0) {\s $1$};
\node[] () at (0.5,0) {\s $2$};
\node[] () at (1,0) {\s $3$};
\draw[<-] (0.1,0) -- (0.4,0);
\draw[<-] (0.6,0) -- (0.9,0);
\draw[->] (0,0.15) arc (120:60:1);
\node[] () at (1.2,0) {$\Bigr)$};
\end{tikzpicture}
\end{center}
\end{minipage}
 and $\mu_1((x_1,\, x_2,\, x_3),Q) =$\break\hspace*{-8pt}
\begin{minipage}{3.88cm}
\begin{center} 
\begin{tikzpicture}
\node[] () at (-1.28,0) {$\Bigl((\frac{1+x_2}{x_3},x_2,x_3),\,$};
\node[] () at (0,0) {\s $1$};
\node[] () at (0.5,0) {\s $2$};
\node[] () at (1,0) {\s $3$};
\draw[<-] (0.1,0) -- (0.4,0);
\draw[->] (0.6,0) -- (0.9,0);
\node[] () at (1.2,0) {$\Bigr)$};
\end{tikzpicture}
\end{center}
\end{minipage}.

We continue by applying $\mu_1,\,\mu_2,\,\mu_3$ to the new seeds, keeping in mind that $\mu_i^2$ is the identity. In this example we get only a finite number of seeds, namely $14.$ Note that the seeds 
\begin{minipage}{3.6cm}
\begin{center} 
\begin{tikzpicture}
\node[] () at (-1.1,0) {$\Bigl((x_1,\, x_2,\, x_3),\,$};
\node[] () at (0,0) {\s $1$};
\node[] () at (0.5,0) {\s $2$};
\node[] () at (1,0) {\s $3$};
\draw[->] (0.1,0) -- (0.4,0);
\draw[->] (0.6,0) -- (0.9,0);
\node[] () at (1.2,0) {$\Bigr)$};
\end{tikzpicture}
\end{center}
\end{minipage}
 and 
\begin{minipage}{3.6cm}
\begin{center} 
\begin{tikzpicture}
\node[] () at (-1.2,0) {$\Bigl((x_1,\, x_3,\, x_2),\,$};
\node[] () at (-0.1,0) {\s $1$};
\node[] () at (0.5,0) {\s $2$};
\node[] () at (1,0) {\s $3$};
\draw[<-] (0.6,0) -- (0.9,0);
\draw[->] (0,-0.1) arc (-130:-60:0.8);
\node[] () at (1.2,0) {$\Bigr)$};
\end{tikzpicture}
\end{center}
\end{minipage} are identified, since they are the same up to relabelling.
\end{example}

\begin{example}\label{expl2}\normalfont
When $Q$ is the quiver 
\begin{minipage}{1.4cm}
\begin{center} 
\begin{tikzpicture}
\node[] () at (1,0) {};
\node[] () at (0,0) {\s $1$};
\node[] () at (1,0) {\s $2$};
\draw[->] (0.1,0) -- (0.9,0);
\end{tikzpicture}
\end{center}
\end{minipage} and $F = \mathbb{Q}(x_1,x_2)$, we have the following complete picture of the graph of seeds, called the \textit{cluster graph.}

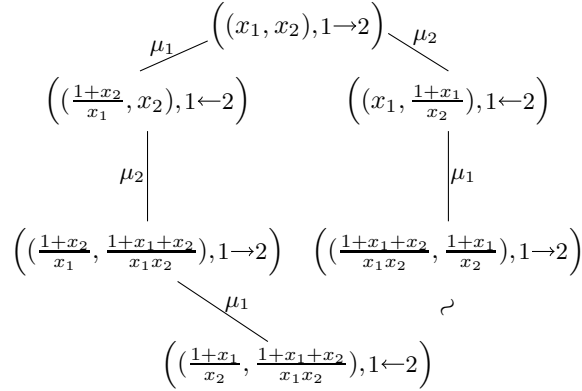
\begin{center}
\begin{tikzpicture}
\node[] (C) at (4,3) 
						{\begin{tikzpicture}
						\node[] () at (-0.8,0) {$\Bigl((x_1,x_2),\,$};
						\node[] () at (0,0) {\s $1$};
						\node[] () at (0.5,0) {\s $2$};
						\draw[->] (0.1,0) -- (0.4,0);
						\node[] () at (0.7,0) {$\Bigr)$};
						\end{tikzpicture}};
\node[] (X) at (2,2) 
						{\begin{tikzpicture}
						\node[] () at (-1.0,0) {$\Bigl((\frac{1+x_2}{x_1},x_2),\,$};
						\node[] () at (0,0) {\s $1$};
						\node[] () at (0.5,0) {\s $2$};
						\draw[<-] (0.1,0) -- (0.4,0);
						\node[] () at (0.7,0) {$\Bigr)$};
						\end{tikzpicture}};
\node[] (B) at (6,2) 
						{\begin{tikzpicture}
						\node[] () at (-1.0,0) {$\Bigl((x_1,\frac{1+x_1}{x_2}),\,$};
						\node[] () at (0,0) {\s $1$};
						\node[] () at (0.5,0) {\s $2$};
						\draw[<-] (0.1,0) -- (0.4,0);
						\node[] () at (0.7,0) {$\Bigr)$};
						\end{tikzpicture}};
\node[] (Y) at (2,0) 
						{\begin{tikzpicture}
						\node[] () at (-1.44,0) 		{$\Bigl((\frac{1+x_2}{x_1},\frac{1+x_1+x_2}{x_1x_2}),\,$};
						\node[] () at (0,0) {\s $1$};
						\node[] () at (0.5,0) {\s $2$};
						\draw[->] (0.1,0) -- (0.4,0);
						\node[] () at (0.7,0) {$\Bigr)$};
						\end{tikzpicture}};
\node[] (A) at (6,0) 
						{\begin{tikzpicture}
						\node[] () at (-1.44,0) 		{$\Bigl((\frac{1+x_1+x_2}{x_1x_2},\frac{1+x_1}{x_2}),\,$};
						\node[] () at (0,0) {\s $1$};
						\node[] () at (0.5,0) {\s $2$};
						\draw[->] (0.1,0) -- (0.4,0);
						\node[] () at (0.7,0) {$\Bigr)$};
						\end{tikzpicture}};
\node[] (Z) at (4,-1.5) 
						{\begin{tikzpicture}
						\node[] () at (-1.44,0) 		{$\Bigl((\frac{1+x_1}{x_2},\frac{1+x_1+x_2}{x_1x_2}),\,$};
						\node[] () at (0,0) {\s $1$};
						\node[] () at (0.5,0) {\s $2$};
						\draw[<-] (0.1,0) -- (0.4,0);
						\node[] () at (0.7,0) {$\Bigr)$};
						\end{tikzpicture}};
\draw[] (2.8,2.8) -- (1.9,2.4);
\node[] () at (2.2,2.7) {\s $\mu_1$};
\draw[] (5.2,2.9) -- (6,2.4);
\node[] () at (5.7,2.85) {\s $\mu_2$};
\draw[] (2,1.6) -- (2,0.4);
\node[] () at (1.8,1) {\s $\mu_2$};
\draw[] (6,1.6) -- (6,0.4);
\node[] () at (6.2,1) {\s $\mu_1$};
\draw[] (2.4,-0.4) -- (3.6,-1.2);
\node[] () at (3.2,-0.7) {\s $\mu_1$};
\node[] () at (6,-0.75) { \rotatebox[origin=c]{45}{$\sim$}};
\end{tikzpicture}
\end{center}
\end{example}

The $n$-element subsets occurring in the seeds are called \emph{clusters,} and the elements occurring in the clusters are called \emph{cluster variables.} Finally, the associated cluster algebra is the subalgebra of $\mathbb{Q}\seqn{x}$ generated by the cluster variables.

In Example \ref{expl2} we have $5$ clusters, and the cluster variables are the $5$ elements $x_1,\, x_2,$   $\frac{1+x_1}{x_2},\,\frac{1+x_2}{x_1},\,\frac{1+x_1+x_2}{x_1x_2}.$ The cluster graph is also interesting in connection with the study of associahedra (\cite{FZ4},\cite{MRZ}).

\subsection{Basic properties}
In both examples above there is only a finite number of seeds, clusters and cluster variables. This is, however, not usually the case. In fact, there is the following nice description of   when this holds \cite{FZ3}.

\begin{theorem}\label{Thm23}
Let $Q$ be a finite connected quiver with no loops or $2$-cycles. Then there are only finitely many clusters (or cluster variables, or seeds) if and only if the underlying graph of $Q$ is a Dynkin diagram.
\end{theorem}

Note that this result is analogous to Gabriel's classification theorem for when there are only finitely many nonisomorphic indecomposable representations of a quiver.

A remarkable property of the cluster variables is the following \cite{FZ1}, called the \textit{Laurent phenomenon.}

\begin{theorem}\label{Thm24}
Let $(\seqn{x},Q)$ be an initial seed.
When we write a cluster variable in reduced form $f/g,$ then $g$ is a monomial in $\seqs{x}.$
\end{theorem}

We shall see later that these monomials $g$ contain some interesting information from a representation theoretic point of view.

A cluster algebra is said to be \textit{acyclic} if there is some seed with an acyclic quiver.  In this case there is the following  information on the numerators of the cluster variables, when they are expressed in terms of the cluster in a seed which has an acyclic quiver (\cite{CK1},\cite{CR},\cite{Q}). The corresponding result is not known for cluster algebras in general.

\begin{theorem}\label{Thm25}
For an acyclic cluster algebra as above, there are positive coefficients for all monomials in the numerator $f$ of a cluster variable in reduced form.
\end{theorem}
\subsection{Cluster algebras with coefficients}\label{clustercoef}
 Cluster algebras with coefficients are important for geometric examples of cluster algebras. Here we only consider a special case of such cluster algebras. Let $\mathbb{Q}(\seqs{x}; \seqs[t]{y})$ be the rational function field in $n+t$ variables, where $\seqs[t]{y}$ are called \textit{coefficients.} Let $Q$ be a finite quiver with $n+t$ vertices corresponding to $\seqs{x}, \seqs[t]{y}.$  We start with the \textit{seed} $((\seqs{x}; \seqs[t]{y}), Q).$ Then we only apply mutations $\seqs{\mu}$ with respect to the first $n$ vertices, and otherwise proceed as before. Note that Theorems \ref{Thm23} and \ref{Thm24} hold also in this setting.

There are several examples of classes of cluster algebras with coefficients, for example the homogeneous coordinate rings of Grassmanians investigated in \cite{Sc} and the coordinate rings $\mathbb{C}[N]$ of unipotent groups (see \cite{GLS1}). There are further examples in (\cite{BFZ},\cite{GLS3},\cite{GLS6},\cite{BIRSc},\cite{GSV1}).

The combinatorics of the cluster algebras in the case of Grassmanians of type $A_n$ can be nicely illustrated geometrically by triangulations of a regular $(n+3)$-gon. The cluster variables, which are coefficients, correspond to the edges of the regular $(n+3)$-gon, and the other cluster variables correspond to the diagonals. The clusters, without coefficients, correspond to the triangulations of the $(n+3)$-gon, or in other words to maximal sets of diagonals which do not intersect. We illustrate with the following simple example.

\begin{example}\normalfont
Let
\begin{minipage}{3cm} 
\begin{center}
\begin{tikzpicture}
\draw[] (0,0.9) -- (1.1,0.2);
\draw[] (1.2,0.2) -- (2,1);
\draw[] (2,1.1) -- (1.5,2);
\draw[] (1.4,2) -- (0.6,2);
\draw[] (0.4,1.9) -- (0,1);
\draw[] (0.1,1) -- (1.4,1.9);
\draw[] (0.1,0.9) -- (1.9,1);
\node[] () at (0.5,2.2) {\s $1$};
\node[] () at (1.5,2.2) {\s $2$};
\node[] () at (-0.2,1) {\s $5$};
\node[] () at (1.1,0) {\s $4$};
\node[] () at (2.2,1) {\s $3$};
\node[] () at (0.7,1.6) {\s $a$};
\node[] () at (1,1.1) {\s $b$};
\end{tikzpicture}
\end{center}
\end{minipage}
 be a regular $5$-gon. Then the diagonals a and b correspond to cluster variables which are not coefficients, and $(a,b)$ is a maximal set of non intersecting diagonals. Taking corresponding cluster variables together with the coefficients, we obtain a cluster. It is easy to see that there are $5$ diagonals, and $5$ triangulations (see \cite{FZ4}).
\end{example}

\section{Cluster Categories}\label{sec3}
It is of interest to categorify the main ingredients in the definition of cluster algebras. The idea behind this is to work within a category with extra structure, and imitate the basic operations within this new category. The hope is on one hand that this will give some feedback to the theory of cluster algebras, and on the other hand that it will lead to new interesting theories. Both aspects have been successful, for various classes of cluster algebras. In this section we deal with the acyclic ones. 

\subsection{Cluster structures}\label{Clusstruc}
We first make a list of desired properties for the categories $\C$ we are looking for. Let $\C$ be a triangulated $k$-category with split idempotents, which is Hom-finite, that is, the homomorphism spaces are finite dimensional over $k.$ Then $\C$ is a Krull-Schmidt category, that is, each object is a finite direct sum of indecomposable objects with local endomorphism ring.  We are looking for  appropriate sets of $n$ nonisomorphic indecomposable objects $\seqs{T},$ or rather objects $T = \dsum{T},$ where the $T_i$ should be the analogs of cluster variables and the $T,$ the analogs of clusters. To have a good analog we would like these objects to satisfy the following, in which case we say that $\C$ has a \textit{cluster structure} (see  \cite{BIRSc}). $\C$ has a \textit{weak cluster structure} if (C1) and (C2) are satisfied.
\begin{itemize}
	\item[(C1)] For $T = \dsum{T}$ in our set, there is, for each $i = \num{n},$ a unique indecomposable object $T_i^*\nsimeq T_i$ in $\C$ such that $T/T_i\oplus T_i^*$ is in our set.
	\item[(C2)] For each $T_i$ there are triangles $T_i^*\ts{f}B_i\ts{g}T_i\rightarrow T_i^*[1]$ and $T_i\ts{s}B_i'\ts{t}T_i^*\rightarrow T_i[1],$ where the maps $g$ and $t$ are minimal right $\add(T/T_i)$-approximations and the maps $f$ and $s$ are minimal left $\add(T/T_i)$-approximations.
	\item[(C3)] There are no loops or $2$-cycles in the quiver $Q_T$ of $\EndopC{T}.$ This means that any nonisomorphism $u:\map{T_i}{T_i}$ factors through $g:\map{B_i}{T_i}$ and through $s:\map{T_i}{B_i'},$ and $B_i$ and $B_i'$ have no common nonzero summands. 
	\item[(C4)] For each $T$ in our set we have $\mu_i(Q_T) = Q_{\mu_i(T)}.$
\end{itemize}

We recall that the map $g:\map{B_i}{T_i}$ is a right $\add(T/T_i)$-approximation if $B_i$  is in $\add(T/T_i)$ and any map $h:\map{X}{T_i}$ with $X$ in $\add(T/T_i)$  factors through $g:\map{B_i}{T_i}.$ The map $g:\map{B_i}{T_i}$  is right minimal if for any commutative diagram  
\begin{minipage}{2cm}
\begin{center} 
\begin{tikzpicture}
\node[] (A) at (0,0) {\s $B_i$};
\node[] (B) at (1.5,0) {\s $T_i$};
\node[] (C) at (0.2,-1) {\s $B_i$};
\draw[->] (0.2,0) -- (1.3,0);
\draw[->] (0,-0.2) -- (0.2,-0.8);
\draw[->] (0.4,-0.8) -- (1.3,-0.2);
\node[] () at (0.75,0.15) {\s $g$};
\node[] () at (1,-0.6) {\s $g$};
\node[] () at (0.2,-0.4) {\s $s$};
\end{tikzpicture}
\end{center}
\end{minipage}, the map $s$ is an isomorphism.  Left approximations and left minimal maps are defined similarly.

Note that the relationship between $T_i$ and $T_i^*$ required in part (C2) is similar to the formula  $x_ix_i^* = m_1+m_2$ appearing in the definition of cluster algebras.

Part (C4) is related to our discussion in Section \ref{sec1} about what is needed in order to categorify quiver mutation. There we saw that tilting $kQ$-modules could not always be used for this purpose. However, for a class of complete algebras of Krull dimension $3,$ known as \textit{$3$-Calabi-Yau algebras,} we shall see that actually the tilting modules (of projective dimension at most $1)$ can be used.

When we deal with cluster algebras with coefficients, we need a modified version of the above definition of cluster structure (see \cite{BIRSc}).

\subsection{Origins}
When $Q$ is a Dynkin quiver, it was shown in \cite{FZ1} that there is a one-one correspondence between the cluster variables for the cluster algebra determined by the seed $(\seqn{x},Q),$ and the almost positive roots, that is, the positive roots together with the negative simple roots. On the other hand the positive roots are in $1-1$ correspondence with the indecomposable $kQ$-modules. This led the authors of \cite{MRZ} to introduce the category of \textit{decorated} representations of $Q,$ which are the representations of the quiver 
\begin{minipage}{2.4cm}
\begin{center} 
\begin{tikzpicture}
\node[] () at (0,0.1) {$Q\cup\bigl\{$};
\node[] () at (0.5,0.2) {\s $1$};
\node[dot] () at (0.5,0) {};
\node[] () at (1,0) {$\ldots$};
\node[] () at (1.5,0.2) {\s $n$};
\node[dot] () at (1.5,0) {};
\node[] () at (1.7,0.1) {$\bigr\},$};
\end{tikzpicture}
\end{center}
\end{minipage} that is, $Q,$ together with $n$ isolated vertices.
Then the indecomposable representations of $Q$ correspond to the
positive roots and the $n$ additional $1$-dimensional representations correspond to the negative simple roots. A compatibility degree $E(X,Y)$ for a pair of decorated representations was introduced, and $E(X,X) = 0$ corresponds to $\ext_{kQ}^1(X,X)=0$ when $X$ is in $\m kQ.$ The maximal such $X$ in $\m kQ$ are the tilting modules, indicating a connection with tilting theory.  On the other hand, if the cluster variables should correspond to indecomposable objects in the category we are looking for, then we would need additional indecomposable objects compared to $\m kQ.$ This could also help remedy the fact that not all almost complete tilting modules in $\m kQ$ have exactly two complements in $\m kQ.$ Recall that the \textit{almost complete} tilting modules are the modules obtained from tilting modules $T = \dsum{T},$ where $n$ is the number of vertices in the quiver,  by dropping one indecomposable summand. In order to have property (C4), we would  also need more maps in our desired category.

Taking all this into account, it turns out to be fruitful to consider a suitable orbit category of the bounded derived category $\mathcal{D}^b(kQ).$ In the case of path algebras the bounded derived categories have a particularly nice structure.

\subsection{Derived categories of path algebras}
Let $Q$ be a finite acyclic quiver, for example 
\begin{minipage}{1.4cm}
\begin{center} 
\begin{tikzpicture}
\node[] () at (0,0) {\s $1$};
\node[] () at (0.5,0) {\s $2$};
\node[] () at (1,0) {\s $3$};
\draw[->] (0.1,0) -- (0.4,0);
\draw[->] (0.6,0) -- (0.9,0);
\end{tikzpicture}
\end{center}
\end{minipage} as before. Then the indecomposable objects in the bounded derived category $\mathcal{D}^b(kQ)$ are just the objects $X[i]$ for $i\in\mathbb{Z},$ where $X$ is an indecomposable $kQ$-module. Then for $X$ and $Y$ in $\mkq$ we have that $\Hom_{\Dbk}(X[i],Y[j])$ is isomorphic to $\Hom_{kQ}(X,Y)$ if $i = j,$ to $ \ext_{kQ}^1(X,Y)$ if $j=i+1,$ and is $0$ otherwise.

The category $\Dbk$ is a triangulated category, and it has \textit{almost split triangles} \cite{Ha2}, where those inside a given shift $(\mkq)[i]$ are induced by the almost split sequences in $\mkq.$ The others are of the form $I[j-1]\rightarrow E\rightarrow P[j],$ where $P$ is indecomposable projective in $\mkq,$ and $I$ is the indecomposable injective module associated with the same vertex of the quiver. The operation $\tau$ is defined as for almost split sequences. Actually, in this setting it is even induced by an equivalence of categories $\tau:\map{\Dbk}{\Dbk}$ \cite{Ha2}. 

For the running example we then have the corresponding AR-quiver
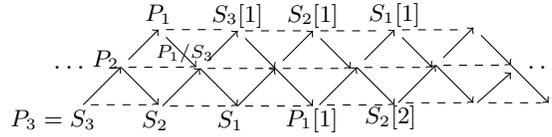
\begin{center}
\begin{tikzpicture}

\node[] () at  (-1,0) {$\ldots$};

\node[] () at  (-0.4,0)
{\begin{tikzpicture}
\draw[->] (0,0) -- (0.5,0.5);
\draw[->] (0.6,0.6) -- (1,1);
\draw[->] (0.5,0.5) -- (1,0);%
\draw[->] (1.1,0.9) -- (1.5,0.5);
\draw[dashed] (1.1,1) -- (1.9,1);
\draw[dashed] (0.6,0.5) -- (1.4,0.5);
\draw[dashed] (0.1,0) -- (0.9,0);
\node[] () at  (-0.4,-0.2) {\s $P_3=S_3$};
\node[] () at  (1,1.2) {\s $P_1$};
\node[] () at  (0.3,0.6) {\s $P_2$};
\end{tikzpicture}};

\node[] () at  (1,0)
{\begin{tikzpicture}
\draw[->] (0,0) -- (0.5,0.5);
\draw[->] (0.6,0.6) -- (1,1);
\draw[->] (0.5,0.5) -- (1,0);
\draw[->] (1.1,0.9) -- (1.5,0.5);
\draw[dashed] (1.1,1) -- (1.9,1);
\draw[dashed] (0.6,0.5) -- (1.4,0.5);
\draw[dashed] (0.1,0) -- (0.9,0);
\node[] () at  (1,1.2) {\s $S_3[1]$};
\node[] () at  (-0.1,-0.2) {\s $S_2$};
\node[] () at  (0.3,0.7) {\scriptsize $P_1/S_3$};
\end{tikzpicture}};

\node[] () at  (2,0)
{\begin{tikzpicture}
\draw[->] (0,0) -- (0.5,0.5);
\draw[->] (0.6,0.6) -- (1,1);
\draw[->] (0.5,0.5) -- (1,0);
\draw[->] (1.1,0.9) -- (1.5,0.5);
\draw[dashed] (1.1,1) -- (1.9,1);
\draw[dashed] (0.6,0.5) -- (1.4,0.5);
\draw[dashed] (0.1,0) -- (0.9,0);
\node[] () at  (1,1.2) {\s $S_2[1]$};
\node[] () at  (-0.1,-0.2) {\s $S_1$};
\end{tikzpicture}};

\node[] () at  (3,0)
{\begin{tikzpicture}
\draw[->] (0,0) -- (0.5,0.5);
\draw[->] (0.6,0.6) -- (1,1);
\draw[->] (0.5,0.5) -- (1,0);
\draw[->] (1.1,0.9) -- (1.5,0.5);
\draw[dashed] (1.1,1) -- (1.9,1);
\draw[dashed] (0.6,0.5) -- (1.4,0.5);
\draw[dashed] (0.1,0) -- (0.9,0);
\node[] () at  (1,1.2) {\s $S_1[1]$};
\node[] () at  (-0.1,-0.2) {\s $P_1[1]$};
\end{tikzpicture}};

\node[] () at  (4.1,0)
{\begin{tikzpicture}
\draw[->] (0,0) -- (0.5,0.5);
\draw[->] (0.6,0.6) -- (1,1);
\draw[->] (0.5,0.5) -- (1,0);
\draw[->] (1.1,0.9) -- (1.5,0.5);
\draw[->] (1.05,0) -- (1.5,0.4);
\draw[->] (1.5,0.5) -- (2,0);
\draw[dashed] (0.6,0.5) -- (1.4,0.5);
\draw[dashed] (0.1,0) -- (0.9,0);
\draw[dashed] (1.2,0) -- (1.9,0);
\node[] () at  (1,1.3) {};
\node[] () at  (-0.1,-0.2) {\s $S_2[2]$};
\end{tikzpicture}};

\node[] () at  (5.3,0) {$\ldots$};

\end{tikzpicture}
\end{center}

When $Q'$ is obtained from $Q$ by reflection at a sink (for example\break\hspace*{-6.5pt}
\begin{minipage}{2.4cm}
\begin{center} 
\begin{tikzpicture}
\node[] () at (0,0) {\s $1$};
\node[] () at (1,0) {\s $2$};
\node[] () at (2,0) {\s $3$};
\draw[->] (0.1,0) -- (0.9,0);
\draw[<-] (1.1,0) -- (1.9,0);
\end{tikzpicture}
\end{center}
\end{minipage} is obtained from 
\begin{minipage}{2.4cm}
\begin{center} 
\begin{tikzpicture}
\node[] () at (0,0) {\s $1$};
\node[] () at (1,0) {\s $2$};
\node[] () at (2,0) {\s $3$};
\draw[->] (0.1,0) -- (0.9,0);
\draw[->] (1.1,0) -- (1.9,0);
\end{tikzpicture}
\end{center}
\end{minipage} in the above example), then we have seen that we have an associated tilting $kQ$-module $T$ such that $\Endop[kQ]{T}\simeq kQ'.$ Hence there is induced a derived equivalence between $kQ$ and $kQ'$ \cite{Ha2}.

\subsection{The cluster category}

%Denote by $[1]:\map{\Dbk}{\Dbk}$ the shift functor. 

Let as before $\tau:\map{\Dbk}{\Dbk}$ denote the  equivalence which
induces the AR-translation, so that $\tau C = A$ if $A\rightarrow
B\rightarrow C\rightarrow A[1]$ is an almost split triangle.   
Then $F = \tau^{-1}[1]:\map{\Dbk}{\Dbk}$ is an equivalence, where [1]
denotes the shift functor. We defined in \cite{BMRRT} the cluster category $\CQ$ to be the orbit category $\Dbk/F$. By definition, the indecomposable objects are the $F$-orbits in $\Dbk$ of indecomposable objects, represented by indecomposable objects in the \textit{fundamental domain} $\mathcal{F}.$ This consists of the indecomposable objects in $\mkq,$ together with $P_1[1],\ldots,P_n[1],$ where $\seqs{P}$ are the indecomposable projective $kQ$-modules. So in the Dynkin case, the number of indecomposable objects in $\mathcal{F},$ up to isomorphism, equals the number of cluster variables. For $X$ and $Y$ indecomposable in $\mathcal{F}$ we have $\Hom_{\CQ}(X,Y) = \underset{i=-\infty}{\overset{\infty}{\bigoplus}}\Hom_{\Db}(X,F^iY),$ by the definition of maps in an orbit category. So in general we have more maps than before. For example, if $Q$ is 
\begin{minipage}{1.4cm}
\begin{center} 
\begin{tikzpicture}
\node[] () at (0,0) {\s $1$};
\node[] () at (0.5,0) {\s $2$};
\node[] () at (1,0) {\s $3$};
\draw[->] (0.1,0) -- (0.4,0);
\draw[->] (0.6,0) -- (0.9,0);
\end{tikzpicture}
\end{center}
\end{minipage}, we have $\Hom_{\Db}(S_1,F(S_3)) \simeq k,$ and
$\Hom_{\Db}(S_1,F^iS_3) = 0$ for $i\neq 1.$  Hence we have
$\Hom_{\CQ}(S_1,S_3)\simeq k.$ Now, for $T = P_1\oplus S_1\oplus S_3$ in $\CQ,$ we have $Q_T =:$ 
\begin{minipage}{2.3cm}
\begin{center} 
\begin{tikzpicture}
\node[] () at (0,0) {\s $1$};
\node[] () at (1,0) {\s $2$};
\node[] () at (2,0) {\s $3$};
\draw[<-] (0.1,0) -- (0.9,0);
\draw[<-] (1.1,0) -- (1.9,0);
\draw[->] (0.1,-0.1) arc (-120:-60:1.8);
\end{tikzpicture}
\end{center}
\end{minipage}, which coincides with $\mu_2(Q),$ as desired. 

The category $\CQ$ has some nice additional structure, which orbit categories rarely have; namely they are triangulated categories \cite{K1}. 

\subsection{Cluster tilting objects}\label{sec35}
In Section \ref{genquivmut} we have seen that we could not use tilting $kQ$-modules to categorify quiver mutation at vertex $2$ in our running example, since the quiver of the tilting module $T = P_1\oplus S_1\oplus P_3$ was not ``correct.'' But when $T$ is viewed as an object in $\CQ,$ the associated quiver is the correct one. So it seems natural that our desired class of objects should include the tilting $kQ$-modules. Now, as already pointed out, there are usually different module categories of the form $\mkq$ giving rise to the same bounded derived category $\mathcal{D}^b(kQ),$ and hence to the same cluster category. So it is natural to consider the objects in $\CQ$ induced by all the corresponding tilting modules. This class of modules turns out to have a nice  uniform description  as objects in $\CQ,$ which motivates the following definitions.

An object $T$ in $\CQ$ is \textit{maximal rigid} if $\ext_{\CQ}^1(T,T)=0$ and $T$ is maximal with this property. It is \textit{cluster tilting} if it is rigid and $\ext_\C^1(T,X) = 0$ implies that $X$ is in $\add T.$ Actually the concepts of maximal rigid and cluster tilting coincide for $\CQ$ \cite{BMRRT}.

We have the following (\cite{BMRRT},\cite{BMR2}).

\begin{theorem}
The cluster category $\CQ$ has  a cluster structure with respect to the cluster tilting objects.
\end{theorem}

As for the cluster algebras, there is also in this setting a natural associated graph, called the \textit{cluster tilting graph.} Associated with a cluster tilting object $T$ in $\CQ,$ we have a quiver $Q_T,$ which is the quiver of $\EndopCQ{T},$ and hence we have a natural \textit{tilting seed} $(T,Q_T)$ \cite{BMR2}.  The vertices of the cluster tilting graph correspond to the cluster tilting objects up to isomorphism, or equivalently, to the tilting seeds. We have an edge between two vertices if the two corresponding cluster tilting objects differ by exactly one indecomposable summand. 

Whereas the cluster graph by definition is connected, this is not automatic for the cluster tilting graph. However, it can be shown to be the case for cluster categories (\cite{BMRRT}, using \cite{HU2}), see also \cite{Hub2}.

\subsection{Categorification of quiver mutation}\label{Catqmut}
Note that we have in particular obtained a way of categorifying
quiver mutation beyond the case of mutation at a sink as discussed in Section \ref{sec1}. So we isolate the more general statement as follows.

\begin{theorem}
Let $Q$ be a finite acyclic quiver, and $Q'$ a quiver obtained from $Q$ by a finite sequence of mutations. Let $i$ be a vertex of $Q'.$ Then there is a cluster tilting object $T'$ in $\CQ$ such that for $T'' = \mu_i(T')$ we have the commutative diagram
\begin{center}
\begin{tikzpicture}
\node[] () at (0,0) {\s $T'$};
\node[] () at (3,0) {\s $T''$};
\node[] () at (-0.4,-1) {\s $Q_{T'}=Q'$};
\node[] () at (2.65,-1) {\s $\mu_i(Q')=Q_{T''}$};
\draw[->] (0.2,0) -- (2.8,0);
\draw[-] (0.2,-0.1) -- (0.2,0.1);
\draw[->] (0.4,-1) -- (1.6,-1);
\draw[-] (0.4,-1.1) -- (0.4,-0.9);
\draw[->] (0,-0.2) -- (0,-0.8);
\draw[-] (-0.1,-0.2) -- (0.1,-0.2);
\draw[->] (2.9,-0.2) -- (2.9,-0.8);
\draw[-] (3,-0.2) -- (2.8,-0.2);
\node[] () at (1.5,0.2) {\s $\mu_i$};
\node[] () at (1,-0.8) {\s $\mu_i$};
\end{tikzpicture}
\end{center}
\end{theorem}

As an illustration of how such a categorification can be useful we state the following result \cite{BR2}.

\begin{theorem}
Let $Q$ be a finite acyclic quiver. Then the mutation class of $Q$ is finite if and only if $Q$ has at most two vertices, or $Q$ is a Dynkin or an extended Dynkin diagram.
\end{theorem}

The point of the categorification is that since the cluster tilting objects are closely related to the tilting $kQ$-modules, we can take advantage of the well developed theory of tilting modules over finite dimensional algebras. The problem amounts to deciding when only a finite number of quivers occur as quivers associated with tilting modules.

We point out that a classification of finite mutation type has recently been obtained in general for finite quivers without loops or $2$-cycles \cite{FST}.

As indicated before, quiver mutation can be categorified using
tilting modules of projective dimension at most $1$ for  
a class of algebras of Krull 
dimension 3 called 3-Calabi Yau algebras (see \cite{IR1})

%$3$-Calabi-Yau algebras \cite{IR1}.

\begin{theorem}
Let $\Lambda$ be a basic $3$-Calabi-Yau algebra given by a quiver $Q$ with relations. Assume that  $Q$ has no loops or $2$-cycles. Then $\mu_i(Q)$ is obtained from the quiver of $\mu_i(\Lambda)$ by removing all $2$-cycles, with $\mu_i$ as defined below.
\end{theorem}

Let $\Lambda = \dsum{P},$ where the $P_i$ are indecomposable projective $\Lambda$-modules, and $P_r\nsimeq P_s$ for $r\neq s$ since $\Lambda$ is basic. There is for a given $i = \num{n}$ a unique indecomposable $\Lambda$-module $P_i^*$ such that $\Lambda/P_i\oplus{P_i^*}$ is a tilting module of projective dimension at most $1,$ and we let $\mu_i(\Lambda) = \Lambda/P_i\oplus{P_i^*}.$

\subsection{A geometric description}
When  $Q$ is a quiver of type $A_n,$ there is an independent categorification of the corresponding cluster algebra along very different lines \cite{CCS}. This is based upon the example discussed in Section~\ref{clustercoef}. We consider the triangulations of the regular $(n+3)$-gon, without including the coefficients, which correspond to the edges in the $(n+3)$-gon. A category with indecomposable objects corresponding to the diagonals in the $(n+3)$-gon was defined in \cite{CCS}. The authors showed that this category is equivalent to the cluster category of type $A_n.$ So we get an interesting geometric description of the cluster category in the $A_n$ case. There is also further work in this direction for $D_n$ \cite{Sci}.

Cluster structures in the context of Teichm\"uller spaces were discussed in (\cite{FG1},\cite{FG3},\cite{FG2}, \cite{GSV1},\cite{GSV2}). This inspired the systematic study of cluster algebras coming from oriented Riemann surfaces with boundary and marked points (\cite{FST},\cite{FT}). Also in this case clusters are in bijection with triangulations. It is easy to see that the mutation class is always finite for these examples.
 
\subsection{Hereditary categories}\label{sec38}
The theory of cluster categories also works when we replace $\mkq$ by an arbitrary Hom-finite hereditary abelian category with tilting object \cite{HRS}, as pointed out in \cite{BMRRT}. It has been shown in \cite{BKL} that in the \textit{tubular} case the cluster tilting graph is connected.

\subsection{$(m)$-cluster categories}\label{sec39}
There is a natural generalization of the cluster categories $\CC{Q} = \Dbk/\tau^{-1}[1]$ to $(m)$-cluster categories $\CC{Q}^{(m)} =  \Dbk/\tau^{-1}[m],$ for $m\geq 1.$ Then $\CC{Q}^{(m)}$ is Hom-finite and also triangulated \cite{K1}. Some more results on cluster categories remain true in the more general setting.

We recall some work from (\cite{T},\cite{W},\cite{Z2},\cite{ZZ}). The concepts of maximal rigid and cluster tilting have a natural generalization to $(m)$-maximal rigid and $m$-cluster tilting objects in $\CC{Q}^{(m)}.$ Also in this setting the concepts coincide. Further, the number of nonisomorphic indecomposable summands of an $m$-cluster-tilting object equals the number of vertices in the quiver $Q.$ If we drop one indecomposable summand from an $m$-cluster tilting object $T,$ there are exactly $m$ different ways to replace it by an indecomposable object, such that we still have an $m$-cluster tilting object.

It was shown in \cite{BT} that also for arbitrary $m$ there is a combinatorial description of mutation of $m$-cluster tilting objects in $m$-cluster categories. In this connection the concept of \textit{coloured} quiver mutation is introduced. There is a geometric description of the $m$-cluster categories for quivers of type $A_n$ and $D_n$ (see \cite{BaMa}). In the Dynkin case the concept of \textit{$m$-clusters} has been introduced in \cite{FR}, and it was shown in (\cite{T},\cite{Z2}) that the $m$-cluster category provides a categorification. This was used in (\cite{T},\cite{Z2}) to simplify proofs of results about $m$-clusters in \cite{FR}.

\section{Calabi-Yau Categories of Dimension Two}\label{sec4}
A crucial property for the investigation of cluster tilting objects in cluster categories $\CQ$ was the functorial isomorphism $D\ext_{\CQ}^1(A,B)\simeq \ext_{\CQ}^1(B,A),$ where $D = \Hom_k(\,,k),$ which by definition expresses that the Hom-finite triangulated $k$-category $\CQ$ is $2$-Calabi-Yau (\CY for short). A similar theory worked for $\m\Lambda$ when $\Lambda$ is the preprojective algebra of a Dynkin diagram \cite{GLS1}.  Also in this case an important feature was that the stable category $\underline{\m}\Lambda$ is \CY.  Then  we say that $\m\Lambda$ is \textit{stably \mbox{\CY.}} This motivated trying to generalize work from cluster categories to arbitrary Hom-finite triangulated \CY $k$-categories. We usually omit the $k$ when we speak about $k$-categories.

\subsection{Preprojective algebras of Dynkin type}
In their work in \cite{GLS1} on $\m\Lambda$ for $\Lambda$ a preprojective algebra of Dynkin type, Geiss-Leclerc-Schr\"oer dealt with the maximal rigid modules, as defined in Section \ref{sec3}. One can here go back and forth between exact sequences in $\m\Lambda$ and triangles in $\underline{\m}\Lambda,$ so for the general theory one can deal with either one of these categories. They followed the same basic outline as in (\cite{BMRRT},\cite{BMR2}), and proved that $\underline{\m}\Lambda$ has a cluster structure in the terminology of Section \ref{sec3}. For (C2) the same proof as for cluster categories could be used, but in the other cases new proofs were necessary. This was also the case for showing that the concepts of maximal rigid and cluster tilting coincide also in this context. For this work some of the results of Iyama on  a higher theory of almost split sequences and Auslander algebras were useful (\cite{I1},\cite{I2}). Actually, in this work Iyama introduced independently the concept of maximal $1$-orthogonal, which coincides with cluster tilting in the setting of (stably) \CY categories. 

The cluster algebras $\mathbb{C}[N]$ are categorified using the cluster tilting objects in the category $\m\Lambda$ for a preprojective algebra $\Lambda$ of Dynkin type. All the indecomposable projective modules are summands of any cluster tilting object, and correspond to the coefficients of the associated cluster algebra. Actually, here categorification can be used to show that $\mathbb{C}[N]$ has a cluster algebra structure (see \cite{GLS1}). This has recently been generalized in \cite{GLS6}.

Cluster monomials are monomials of cluster variables in a given cluster. A central question is their relationship to the canonical and semicanonical bases and their duals, investigated by Lusztig and Kashiwara (see \cite{GLS1},\cite{GLS5}).

\subsection{Generalizations}
In the general case of stably \CY categories, or Hom-finite triangulated \CY categories, there are not necessarily any cluster tilting objects. Actually, there may be maximal rigid objects, but no cluster tilting objects \cite{BIKR}. But we have the following general result \cite{IY}.

\begin{theorem}
Let $\C$ be a $\Hom$-finite triangulated \CY category with cluster tilting objects. Then $\C$ has a weak cluster structure.
\end{theorem}
The proof of (C2) is the same as for cluster categories, while a new argument was needed for (C1). Property (C3) does not however hold in general. There are many examples of stable categories of Cohen-Macaulay modules over isolated hypersurface singularities where there are both loops and $2$-cycles \cite{BIKR}. But if there are no loops or $2$-cycles, then we have the following \cite{BIRSc}.

\begin{theorem}
Let $\C$ be a $\Hom$-finite \CY triangulated $k$-category having cluster tilting objects, and with no loops or 2-cycles. Then $\C$ has a cluster structure.
\end{theorem}

We have pointed out that the cluster tilting graph is known to be connected for cluster categories $\CQ$ when $Q$ is a finite connected quiver. This is an important open problem for connected Hom-finite triangulated \CY categories in general. The only other cases where this is known to be true is for the case discussed in Section \ref{sec38}, and for some cases of stable categories of Cohen-Macaulay modules in \cite{BIKR}.

\subsection{\CY categories associated with elements in Coxeter groups}
Let $Q$ be a finite acyclic quiver  with vertices $\num{n},$ and let $C_Q$ be the associated Coxeter group. By definition $C_Q$ has a distinguished set of generators $\seqs{s},$ and the relations are as follows: $s_i^2 = 1,\,\, s_is_j = s_js_i$ if there is no arrow between $i$ and  $j$ in $Q,$ and $s_is_js_i = s_js_is_j$ if there is exactly one arrow between $i$ and $j.$  Let $\mathbf{w} = s_{i_1}\ldots s_{i_t}$ be a \textit{reduced expression} of $w$ in $\CQ,$ that is, $t$ is smallest possible. Let $\Lambda$ be the (completion of the) preprojective algebra associated with $Q.$ For each $i = \num{n},$ consider the ideal $I_i = \Lambda(1-e_i)\Lambda$ in $\Lambda,$ where $e_i$ denotes the trivial path at the vertex $i.$ Then define the ideal $I_w= I_{i_1}\ldots I_{i_t}.$ It can be shown to be independent of the reduced expression, and the factor algebra $\Lambda_w = \Lambda/I_w$ is finite dimensional. Denote by $\Sub$ the full subcategory of $\m\Lambda_w$ whose objects are the submodules of the free $\Lambda_w$-modules of finite rank. Then $\C_w = \Sub$ is stably \CY, and the associated stable category $\Subu$ is Hom-finite triangulated \CY. Here we refer to (\cite{IR1},\cite{BIRSc}). 

There is a nice class of cluster tilting objects, called \textit{standard} cluster tilting objects, in $\Sub$ and $\Subu.$ When $\mathbf{w} = s_{i_1}\ldots s_{i_t}$ is a reduced expression we define 
$$T_{\mathbf{w}} = \left(P_{i_1}/I_{i_1}P_{i_1}\right)\oplus \left(P_{i_2}/I_{i_1}I_{i_2}P_{i_2}\right)\oplus\ldots\oplus \left(P_{i_t}/I_{i_1}\ldots{}I_{i_t}P_{i_t}\right).$$
Then $T_\mathbf{w}$ is a cluster tilting object in $\Sub$ and in $\Subu,$ and it depends on the reduced expression. However, the standard cluster tilting objects all lie in the same component of the cluster tilting graph \cite{BIRSc}.

A stably \CY category, dual to $\C_w,$ was independently associated with a class of words called \textit{adaptable,} in a very different way \cite{GLS4}. Here there were also associated two cluster tilting objects in a natural way, which are a subset of the set of standard cluster tilting objects (up to duality).

\subsection{Stable categories of Cohen-Macaulay modules}
Let $R$ be a commutative  complete local Gorenstein isolated singularity of Krull dimension $3,$ with $k\subset R$ for the field $k.$ Then by  results of Auslander \cite{Au} on existence of almost split sequences for (maximal) Cohen-Macaulay modules, the category of  Cohen-Macaulay modules $\CM(R)$ is stably \CY and the stable category $\underline{\CM}(R)$ is Hom-finite triangulated \CY.

When $R$ is an isolated hypersurface singularity, we can, by the periodicity result for hypersurfaces (\cite{Kn},\cite{So}),  deal with the case of Krull dimension $1$ just as well.  Already for finite Cohen-Macaulay type there are examples with no cluster tilting objects, and where there are maximal rigid objects which are not cluster tilting \cite{BIKR}. For all these examples there are $2$-cycles, and in many cases also loops. An interesting question in this connection is the following.
Let $\C$ be a Hom-finite triangulated \CY $k$-category, where we have no loops or $2$-cycles. Then do the maximal rigid objects coincide with the cluster tilting objects?

Another class of Gorenstein rings giving rise to \CY categories with desired properties is the following. Let $G\subset SL(3,k)$ be a finite subgroup where no $g\neq 1$ in $G$ has eigenvalue $1,$ and let $R = k[[X,Y,Z]]^G$ be the associated invariant ring, which under these assumptions is an isolated singularity. Then $\underline{\CM}(R)$ is Hom-finite triangulated \CY, and does not have loops or $2$-cycles \cite{KNi}.

\subsection{Generalized cluster categories}
In \cite{Am1} Amiot introduced a new class of Hom-finite triangulated \CY-categories, generalizing the class of cluster categories.

Let $A$ be a finite dimensional $k$-algebra of global dimension at most $2.$ Also under this assumption $\D{A}$ has almost split triangles, and the AR-translation $\tau$ is induced by an equivalence $\tau:\map{\D{A}}{\D{A}}$ \cite{Ha2}.  Consider again the orbit category $\D{A}/\tau^{-1}[1].$  In this setting the orbit category is not necessarily triangulated. The \textit{generalized cluster category} $\CC{A}$ is then defined to be the \textit{triangulated hull} of $\D{A}/\tau^{-1}[1].$
If $\Hom_{\CC{A}}(A,A)$ is finite dimensional, then $\CC{A}$ is Hom-finite triangulated \CY, with $A$ as a cluster tilting object (see \cite{Am1}).

A striking application of these generalized cluster categories is Keller's proof of the periodicity conjecture for pairs of Dynkin diagrams \cite{K5} (see also \cite{IIKNS}), using the algebra $A = kQ\otimes_kkQ$ of global dimension at most $2,$ where $Q$ is a Dynkin quiver.

A more general construction of Hom-finite triangulated \CY categories was given in (\cite{Am1},\cite{K4}), starting with a quiver with potential $(Q,W)$ such that the associated Jacobian algebra is finite dimensional (see Section \ref{sec5}). There is an associated differential graded algebra $\Gamma,$ called the Ginzburg algebra, and a triangulated \CY category $\C_{(Q,W)}$ was constructed from $\Gamma.$ By (\cite{Am1},\cite{K4}) this generalizes the previous construction of $\CC{A}$ from $A.$

\subsection{Relationship between the different classes}
A natural question to ask is how the various classes of Hom-finite triangulated \CY categories are related.

We first show how the cluster categories and the preprojective algebras of Dynkin type are related to the \CY categories associated with elements in Coxeter groups (\cite{BIRSc},\cite{GLS4}).

\begin{theorem}\label{Thm43}
\begin{enumerate}
\item The cluster category $\CQ$ is triangle equivalent to the category $\Subu,$ where $w = c^2$ for a Coxeter element $c,$ when $Q$ is not of type $A_n.$
\item When $Q$ is Dynkin and $\Lambda$ is the associated preprojective algebra, then $\m\Lambda$ is $\Sub,$ where $w$ is the longest element in the Coxeter group.
\end{enumerate}
\end{theorem}
 
The following was shown in (\cite{Am1},\cite{AIRT},\cite{ART}).

\begin{theorem}
Let $\C_w$ be a $\Hom$-finite triangulated \CY category associated with an element in a Coxeter group.  Then there is some algebra $A$ of global dimension at most $2$ such that $\C_w$ is triangle equivalent to the generalized cluster category $\CC{A}.$
\end{theorem}

Another result of a similar flavour was recently shown in (\cite{AIR},\cite{DV}).

\begin{theorem}
Let $R = k[[X,Y,Z]]^G$ be an invariant ring as discussed above, where $G\subset SL(3,k)$ is a finite cyclic group. Then the \CY triangulated category $\underline{\CM}(R)$ is triangle equivalent to a generalized cluster category $\CC{A}$ for some finite dimensional algebra $A$ of global dimension at most $2.$
\end{theorem}

\subsection{Subfactor constructions}
There is also a useful way of constructing new Hom-finite triangulated \CY categories from old ones, via subfactor constructions \cite{IY}.

Let $\C$ be a Hom-finite triangulated \CY category with a nonzero rigid object $D.$ Let ${}^\bot D[1] = \{X\in\C; \ext_\C^1(X,D)=0\}.$ Then the factor category ${}^\bot{D}[1]/\add{D}$ is triangulated \CY. The cluster tilting objects in ${}^\bot{D}[1]/\add{D}$ are in one-one correspondence with the cluster tilting objects of $\C$ which have $D$ as a summand.

This was proved in the general case in \cite{IY}. Related results for the cluster categories were first proved in \cite{BMR2}. There they were used to show that if an algebra $\Gamma$ is cluster tilted, then $\Gamma/\Gamma{e}\Gamma$ is cluster tilted for any idempotent element $e$ of $\Gamma$ (see Section \ref{sec5} for definition). This was useful for reducing problems to algebras with fewer simple modules.

When $\mathbf{w} = \mathbf{uv}$ is a reduced expression in a Coxeter group, then the \CY triangulated category $\underline{\mbox{Sub}}\Lambda_v$ is triangle equivalent to a specific subfactor category of $\Subu$ \cite{IR2}. This was used to get information on components of cluster tilting graphs for $\Subu$ \cite{IR2}.

\section{$2$-Calabi-Yau Tilted Algebras and Jacobian Algebras}\label{sec5}

The study of cluster categories gave rise to an interesting  class of finite dimensional algebras, obtained by taking endomorphism algebras of cluster tilting objects. They have been called \textit{cluster tilted algebras} \cite{BMR1},  and are in some sense analogous to the tilted algebras in the representation theory of finite dimensional algebras.  But their properties are quite different, from several points of view. The cluster tilted algebras have a natural generalization to what has been called \textit{\CY-tilted algebras,} where in the definition we replace cluster categories by Hom-finite triangulated \CY-categories.  In this section we give some basic properties of these algebras, and  discuss their relationship to another important class of algebras; the Jacobian algebras associated with quivers with potential (\cite{DWZ1},\cite{DWZ2}).

\subsection{Special properties of cluster tilted algebras}\label{sec51}
The following result from \cite{ABS} gives a procedure for passing from a tilted algebra to a cluster tilted algebra (see also (\cite{Z},\cite{Ri})). It has no known analog in the general case of \CY-tilted algebras.

\begin{theorem}
Let $\Gamma$ be a tilted algebra. Then the trivial extension algebra $\Gamma\ltimes\ext_\Gamma^2(D\Gamma,\Gamma)$ is cluster tilted.
\end{theorem}

In practice, this is interpreted to amount to drawing an additional arrow from $j$ to $i$ in the quiver of $\Gamma,$ for each relation from $i$ to $j$ in a minimal set of relations for $\Gamma.$ In the case of finite representation type, this construction was first made in \cite{BR1}, and also used in \cite{BRS} for a small class of algebras of infinite representation type. 

For example, if $Q$ is the quiver 
\begin{minipage}{2.4cm}
\begin{center} 
\begin{tikzpicture}
\node[] (A) at (0,0) {\s $1$};
\node[] (B) at (1,0) {\s $2$};
\node[] (C) at (2,0) {\s $3$};
\node[] () at (0.5,0.2) {\s $\a$};
\node[] () at (1.5,0.2) {\s $\b$};
\draw[->] (0.1,0) -- (0.9,0);
\draw[->] (1.1,0) -- (1.9,0);
\end{tikzpicture}
\end{center}
\end{minipage}
and $\Gamma = kQ/\langle\b\a\rangle,$ then the quiver of $\Gamma\ltimes\ext_\Gamma^2(D\Gamma,\Gamma)$ is 
\begin{minipage}{2.4cm}
\begin{center} 
\begin{tikzpicture}
\node[] (A) at (0,0) {\s $1$};
\node[] (B) at (1,0) {\s $2$};
\node[] (C) at (2,0) {\s $3$};
\draw[->] (0.1,0) -- (0.9,0);
\draw[->] (1.1,0) -- (1.9,0);
\draw[<-] (0.1,-0.1) arc (-120:-60:1.8);
\end{tikzpicture}
\end{center}
\end{minipage}

Another interesting property of cluster tilted algebras, not shared by tilted algebras or by \CY-tilted algebras in general, is the following \cite{BIRSm}.

\begin{theorem}\label{thm52}
A cluster tilted algebra is determined by its quiver.
\end{theorem}

In the case of finite representation type this was  proved in \cite{BMR3}, where also the relations were described. Part of the finite type case was also proved independently in \cite{CCS2}.
In general there is no known way of constructing the unique cluster tilted algebra associated with a given quiver.

The following example shows  that the corresponding result does not hold for tilted algebras, since for the quiver $Q:$ 
\begin{minipage}{2.4cm}
\begin{center} 
\begin{tikzpicture}
\node[] () at (0,0) {\s $1$};
\node[] () at (1,0) {\s $2$};
\node[] () at (2,0) {\s $3$};
\draw[->] (0.1,0) -- (0.9,0);
\draw[->] (1.1,0) -- (1.9,0);
\node[] () at (0.5,0.2) {\s $\a$};
\node[] () at (1.5,0.2) {\s $\b$};
\end{tikzpicture}
\end{center}
\end{minipage}, both $kQ$ and $kQ/\langle\b\a\rangle$ are tilted algebras.  If $Q$ is the quiver  
\begin{minipage}{2.4cm}
\begin{center} 
\begin{tikzpicture}
\node[] () at (0,0) {\s $1$};
\node[] () at (1,0.4) {\s $2$};
\node[] () at (2,0) {\s $3$};
\draw[<-] (0.1,0) -- (1.9,0);
\draw[->] (1.1,0.4) -- (1.9,0.1);
\draw[->] (0.1,0.1) -- (0.9,0.4);
\node[] () at (0.4,0.4) {\s $\a$};
\node[] () at (1.6,0.4) {\s $\b$};
\node[] () at (1,-0.2) {\s $\g$};
\end{tikzpicture}
\end{center}
\end{minipage}, then $kQ/\langle\b\a,\g\b,\a\g\rangle$ and $kQ/\langle\b\a\g\b\a,\g\b\a\g\b,\a\g\b\a\g\rangle$ can both be shown to be \CY-tilted algebras, even though they have the same quiver. Here only the first one is cluster tilted. 

\subsection{Homological properties}
The central properties of cluster tilted algebras of a homological nature remain valid also in the general case of \CY-tilted algebras. A Hom-finite triangulated category $\C$ with split idempotents is $3$-CY if we have a functorial isomorphism $D(\ext^1(X,Y))\simeq \ext^2(Y,X)$ for $X,Y$  in $\C.$

\begin{theorem}
Let $\C$ be a $\Hom$-finite triangulated \CY category, $T$ a cluster tilting object in $\C,$ and let $\Gamma = \EndopC{T}.$ Then we have the following.
\begin{enumerate}
	\item The functor $\Hom_\C(T,\,):\map{\C}{\m\,\Gamma}$ induces an equivalence of categories\newline $\C/\add\tau{T}\overset{\sim}{\rightarrow}\m\,\Gamma.$
	\item $\Gamma$ is Gorenstein of injective dimension at most $1.$
	\item The stable category $\underline{\sSub}\,\Gamma$ is triangulated $3$-CY.
\end{enumerate}
\end{theorem}

Part (a) was proved in \cite{BMR1} for cluster tilted algebras and in \cite{KR1} in general. Parts (b) and (c) were proved in \cite{KR1}. 

Part (a) expresses a close relationship between $\C$ and $\m\Gamma.$ For example, on the level  of objects, the indecomposable objects in $\m\Gamma$ are obtained from those in $\C$ by dropping the indecomposable summands of $\tau T$ (which are only finite in number). The category $\C$ has almost split triangles inherited from the almost split triangles in $\Dbk,$ and by dropping the indecomposable summands of $\tau T$ one obtains the AR-quiver of $\m\Gamma.$ Since $\C$ and $\m\Gamma$ are closely related, then also $\m\Gamma$ and $\m\Gamma'$  are closely related when $\Gamma' = \Endop{T'}$ for some cluster tilting object $T'$ in $\C.$ In particular, we have the following (\cite{BMR1},\cite{KR1}), which generalizes Theorem \ref{Thm12}.

\begin{theorem}\label{Thm54}
Let the notation be as above, and assume in addition that $T'$ is obtained from $T$ by a mutation. Then $\Gamma$ and $\Gamma'$ are \emph{nearly} Morita equivalent, that is, there are simple modules  $S$ and $S'$ over $\Gamma$ and $\Gamma'$ respectively, such that the factor categories $\m\Gamma/\add{S}$ and $\m\Gamma/\add{S'}$ are equivalent.
\end{theorem}

Here the objects in $\add{S}$ are finite direct sums of copies of $S,$ and the maps in $\m\Gamma/\add{S}$ are the usual $\Gamma$-homomorphisms modulo those factoring through an object  in $\add S.$ 

In view of the close relationship between $\C$ and $\m\Gamma,$ it is natural to ask if $\m\Gamma$ determines $\C.$ It is not known if this holds in general, but there is the following information \cite{KR2}, which was essential for the proof of Theorem \ref{Thm43}.

\begin{theorem}
Let $\C$ be a $\Hom$-finite triangulated \CY category over the field $k,$ and assume that $\C$ is algebraic (see \cite{KR2}). If there is a cluster tilting object in $\C$ whose associated quiver $Q$ has no oriented cycles, then $\C$ is triangle equivalent to the cluster category $\CQ.$
\end{theorem}

\subsection{Relationship to Jacobian algebras}
It was clear from the beginning of the theory that many examples of cluster tilted algebras, and later of \CY tilted algebras, were given by quivers with potentials. This is a class of algebras appearing in the physics literature \cite{BP}, and they have been systematically investigated in (\cite{DWZ1},\cite{DWZ2}). We refer to \cite{DWZ1} for the general definition of quiver with potential. In particular, a theory of mutations of quivers with potential has been developed.

For example if we have the quiver 
\begin{minipage}{2.4cm}
\begin{center} 
\begin{tikzpicture}
\node[] () at (0,0) {\s $1$};
\node[] () at (1,0.4) {\s $2$};
\node[] () at (2,0) {\s $3$};
\draw[<-] (0.1,0) -- (1.9,0);
\draw[->] (1.1,0.4) -- (1.9,0.1);
\draw[->] (0.1,0.1) -- (0.9,0.4);
\node[] () at (0.4,0.4) {\s $\a$};
\node[] () at (1.6,0.4) {\s $\b$};
\node[] () at (1,-0.2) {\s $\g$};
\end{tikzpicture}
\end{center}
\end{minipage}, we can consider the \textit{potential} $W =\g\b\a,$ which is a cycle. Taking the derivatives with respect to the arrows $\a,\,\b,\,\g,$ up to cyclic permutation, we get $\partial{W}/\partial\a = \g\b,\, \partial{W}/\partial\b = \a\g,\, \partial{W}/\partial \g = \b\a.$  Using these elements as a generating set for the relations, we obtain the first algebra in the  example in Section \ref{sec51}. If we take the potential $W = \g\b\a\g\b\a,$ we get the second algebra. 

The algebras associated with a quiver with potential $(Q,W)$ as above are called \textit{Jacobian} algebras and denoted by $\Jac(Q,W).$ They are not necessarily finite dimensional. For example, the $3$-CY algebras of Krull dimension $3$ mentioned in Section \ref{sec3} are often Jacobian.

The connection between \CY-tilted and Jacobian algebras indicated by the above examples is not accidental. In fact we have the following (\cite{Am1},\cite{K4}).

\begin{theorem}\label{Thm56}
Any finite dimensional Jacobian algebra is \CY-tilted.
\end{theorem}

It is an open problem whether any \CY-tilted algebra is Jacobian. But in many situations this is known to be the case. For example we have the following (\cite{BIRSm},\cite{K4}) (see also \cite{K4} and Corollary \ref{cor512} for more general results).

\begin{theorem}
Any cluster tilted algebra is Jacobian.
\end{theorem}

For the \CY-tilted algebras associated with  standard cluster tilting objects in $\Subu,$ there is an explicit description of the quiver in terms of the reduced expression \cite{BIRSc}.  The same quiver appeared in \cite{BFZ} in the Dynkin case.  In \cite{BIRSm} the following was shown.

\begin{theorem}
Let $w$ be an element in a Coxeter group. Then the \CY-tilted algebras associated with the standard cluster tilting objects in  $\Subu$  are Jacobian.
\end{theorem}

We illustrate with the following.

\begin{example}\normalfont
Let $Q$ be the quiver 
\begin{minipage}{2.5cm}
\begin{center} 
\begin{tikzpicture}
\node[] () at (0,0) {\s $1$};
\node[] () at (1,0) {\s $2$};
\node[] () at (2,0) {\s $3$};
\draw[->] (1.1,0) -- (1.9,0);
\draw[->] (0.1,0) -- (0.9,0);
\node[] () at (0.5,0.2) {\s $a$};
\node[] () at (1.5,0.2) {\s $b$};
\end{tikzpicture}
\end{center}
\end{minipage} and $\mathbf{w} = s_1s_2s_3s_1s_2s_1$ a reduced expression. Let $T$ be the corresponding cluster tilting object in $\Sub.$ Then $\Endop{T}$ and $\underline{\mbox{End}}(T)^{\op}$ have quivers\newline
\begin{minipage}{4.5cm}
\begin{center} 
\hspace*{40pt}\begin{tikzpicture}
\node[] () at (0,0) {\s $1$};
\node[] () at (2,0) {\s $1'$};
\node[] () at (4,0) {\s $1''$};
\node[] () at (1,1) {\s $2$};
\node[] () at (3,1) {\s $2'$};
\node[] () at (2,2) {\s $3$};
\draw[<-] (0.1,0) -- (1.8,0);
\draw[<-] (2.1,0) -- (3.8,0);
\draw[<-] (1.1,1) -- (2.9,1);
\draw[->] (0.1,0.1) -- (0.9,0.9);
\draw[->] (1.1,0.9) -- (1.9,0.1);
\draw[->] (2.1,0.2) -- (2.9,0.9);
\draw[->] (3.1,0.9) -- (3.9,0.1);
\draw[->] (1.1,1.1) -- (1.9,1.9);
\draw[->] (2.1,1.9) -- (2.9,1.1);
\end{tikzpicture}
\end{center}
\end{minipage} \hspace*{40pt}and 
\begin{minipage}{2.5cm}
\begin{center}
\hspace*{30pt}\begin{tikzpicture}
\node[] () at (0,0) {\s $1$};
\node[] () at (2,0) {\s $1'$};
\node[] () at (1,1) {\s $2$};
\draw[<-] (0.1,0) -- (1.8,0);
\draw[->] (0.1,0.1) -- (0.9,0.9);
\draw[->] (1.1,0.9) -- (1.9,0.1);
\node[] () at (0.4,0.6) {\s $a$};
\node[] () at (1.7,0.6) {\s $a^*$};
\node[] () at (1,0.2) {\s $p$};
\end{tikzpicture}
\end{center}
\end{minipage}.

\noindent For the second quiver we have the potential $W = pa^*a,$ and $\underline{\mbox{End}}(T)^{\op}\simeq \Jac(Q,W).$
\end{example}

\subsection{Mutation of quivers with potentials}
Let $(Q,W)$ be a quiver with potential, where $Q$ is a finite quiver with no loops. Then a mutation $\mu_i(Q,W)$ is defined in \cite{DWZ1} (see also \cite{BP}), when $i$ does not lie on a $2$-cycle. Here we illustrate the definition on an example.

\begin{example}\normalfont
Let $Q$ be the quiver 
\begin{minipage}{2.5cm}
\begin{center} 
\begin{tikzpicture}
\node[] () at (0,0) {\s $1$};
\node[] () at (1,0) {\s $2$};
\node[] () at (2,0) {\s $3$};
\draw[->] (1.1,0) -- (1.9,0);
\draw[->] (0.1,0) -- (0.9,0);
\node[] () at (0.5,0.2) {\s $a$};
\node[] () at (1.5,0.2) {\s $b$};
\node[] () at (1.3,-0.4) {\s $c$};
\draw[<-] (0.1,-0.1) arc (-120:-60:1.8);
\end{tikzpicture}
\end{center}
\end{minipage}, and $W = cba$ a potential. We first define $\tilde{\mu}_i(Q,W) = (\tilde{Q},\tilde{W}),$ where $\tilde{Q}$ is the quiver 
\begin{minipage}{2.5cm}
\begin{center} 
\begin{tikzpicture}
\node[] () at (0,0) {\s $1$};
\node[] () at (1,0) {\s $2$};
\node[] () at (2,0) {\s $3$};
\draw[<-] (1.1,0) -- (1.9,0);
\draw[<-] (0.1,0) -- (0.9,0);
\node[] () at (0.5,0.2) {\s $a^*$};
\node[] () at (1.5,0.2) {\s $b^*$};
\node[] () at (1.3,-0.4) {\s $c$};
\node[] () at (1.8,-0.6) {\s $[ba]$};
\draw[<-] (0.1,-0.1) arc (-120:-60:1.8);
\draw[->] (0.1,-0.3) arc (-120:-60:1.8);
\end{tikzpicture}
\end{center}
\end{minipage} and $\tilde{W} = c[ba] + a^*b^*[ba].$ Here we have replaced the path $ba$ of length $2$ going through the vertex $2$ by a new arrow $[ba],$ and we have added a new term $a^*b^*[ba]$ in the potential. Since $\tilde{W}$ has a term which is a cycle of length $2,$ it is by definition not \textit{reduced.} In the next step we get rid of cycles of length $2$ in the potential, and obtain $\mu_2(Q,W) = (\bar{Q},\bar{W}),$ where $\bar{Q}$ is 
\begin{minipage}{2.5cm}
\begin{center} 
\begin{tikzpicture}
\node[] () at (0,0) {\s $1$};
\node[] () at (1,0) {\s $2$};
\node[] () at (2,0) {\s $3$};
\draw[<-] (1.1,0) -- (1.9,0);
\draw[<-] (0.1,0) -- (0.9,0);
\node[] () at (0.5,0.2) {\s $a^*$};
\node[] () at (1.5,0.2) {\s $b^*$};
\end{tikzpicture}
\end{center}
\end{minipage} and $\bar{W} = 0.$
\end{example}

In general $\bar{Q}$ may have $2$-cycles. Then $\bar{Q}$ coincides with $\mu_2(Q)$ only after removing all $2$-cycles. If we require to remove all $2$-cycles, mutation of quivers with potential gives a categorification of quiver mutation.

\subsection{Comparing mutations}
Since there is a large intersection between the classes of Jacobian algebras and \CY-tilted algebras, it makes sense to ask if the mutations  of cluster tilting objects and of quivers with potential are closely related (see \cite{BIRSm}). 

Consider the following diagram, where $T$ is a cluster tilting object in a Hom-finite triangulated \CY category $\C,$ and $(Q,W)$ is a quiver with potential. We assume that $\EndopC{T}\simeq \Jac(Q,W),$ and consider the diagram
\begin{center} 
\begin{tikzpicture}
\draw[->] (0.5,0) -- (2.2,0);
\draw[] (0.5,-0.1) -- (0.5,0.1);
\draw[->] (0.5,2.5) -- (3,2.5);
\draw[] (0.5,2.4) -- (0.5,2.6);
\node[] () at (0,0) {\s $(Q,W)$};
\node[] () at (3.65,0) {\s $\mu_i(Q,W)=(\overline{Q},\overline{W})$};
\node[] () at (0,1) {\s $\Jac(Q,W)$};
\node[] () at (4,1) {\s $\Jac(\overline{Q},\overline{W})$};
\node[] () at (0,1.5) {\s $\Endop{T}$};
\node[] () at (4,1.5) {\s $\Endop{\mu_i(T)}$};
\node[] () at (0,2.5) {\s $T$};
\node[] () at (4,2.5) {\s $\mu_i(T)$};
\draw[->] (0,0.2) -- (0,0.8);
\draw[] (-0.1,0.2) -- (0.1,0.2);
\draw[->] (4,0.2) -- (4,0.8);
\draw[] (3.9,0.2) -- (4.1,0.2);
\draw[->] (0,2.3) -- (0,1.7);
\draw[] (-0.1,2.3) -- (0.1,2.3);
\draw[->] (4,2.3) -- (4,1.7);
\draw[] (3.9,2.3) -- (4.1,2.3);
\node[] () at (0,1.25) { \rotatebox[origin=c]{-90}{$\simeq$}};
\end{tikzpicture}
\end{center}

It is not clear whether any cluster tilting object $T'$ with $\Endop{T'}\simeq\Endop{T}$  gives rise to an algebra $\Endop{\mu_i(T')}$ which is isomorphic to $\Endop{\mu_i(T)}.$ 
Similarly, it is not clear if a quiver with potential $\mu_i(Q',W')$ gives rise to an algebra isomorphic to $\Jac(\mu_i(Q,W))$ when $\Jac(Q,W)\simeq \Jac(Q',W').$ The latter was posed as a problem in \cite{DWZ1}. It was solved in the finite dimensional case, as a consequence of the following \cite{BIRSm}. 

\begin{theorem}\label{Thm511}
Let the notation be as above, and assume that we have an isomorphism $\Endop{T}\simeq\Jac(Q,W).$
\begin{enumerate}
\item For any choice of $T$ and of $(Q,W)$ in the isomorphism $\Endop{T}\simeq\Jac(Q,W),$  we have  $\Endop{\mu_i(T)}\simeq\Jac(\mu_i(Q,W)).$
\item As a consequence, the assignment $\Endop{T}\mapsto\Endop{\mu_i(T)}$ is independent of the choice of $T$ and the assignment 
$\Jac(Q,W)\mapsto\Jac(\overline{Q},\overline{W})$ is independent of the choice of $(Q,W).$
\end{enumerate}
\end{theorem}

We have the following important consequence.

\begin{corollary}\label{cor512}
If a \CY-tilted algebra is Jacobian, then all \CY-tilted algebras belonging to the same component in the cluster tilting graph are Jacobian.
\end{corollary}

Note that this gives an easy proof of the fact that a cluster tilted algebra is Jacobian, since $kQ$ is clearly Jacobian, and we know that there is only one component in this case. 

It also follows, using Corollary \ref{cor512}, that any \CY-tilted algebra belonging to the same component of the cluster tilting graph as those coming from standard cluster tilting objects in categories $\Subu$ are Jacobian. This emphasizes the importance of the problem of proving the existence of only one component in general.

We have seen that for a cluster tilting object $T$ in  a triangulated \CY category $\C,$ then $\Endop{T}$ and $\Endop{\mu_i(T)}$ are nearly Morita equivalent. Hence the corresponding result holds for \textit{neighbouring} finite dimensional Jacobian algebras by Theorems \ref{Thm54}, \ref{Thm56} and \ref{Thm511}. Actually, the following more general result holds (\cite{BIRSm},\cite{DWZ2},\cite{KY}).

\begin{theorem}
If $\Lambda$ and $\Lambda'$ are neighbouring Jacobian algebras, then the categories of finite dimensional modules over $\Lambda$ and $\Lambda'$ are nearly Morita equivalent.
\end{theorem}

The $3$-Calabi-Yau tilted algebras mentioned in Sections \ref{Clusstruc} and \ref{Catqmut} are sometimes given by quivers with potential, and for these algebras we have mutation using tilting modules of projective dimension at most $1.$ Also in this setting the mutation of quivers with potential is closely related to the mutation of tilting modules \cite{BIRSm}.

\begin{theorem}
Assume that $T$ is a tilting module of projective dimension at most $1$ over a $3$-CY-algebra, where $\Endop{T}\simeq\Jac(Q,W).$ Then $\Endop{\mu_i(T)}\simeq\Jac(\mu_i(Q,W)).$
\end{theorem}

In particular, it follows that $\Jac(\mu_i(Q,W))$ is obtained from $\Jac(Q,W)$ via a tilting module $T$ of projective dimension at most $1.$

\subsection{Derived equivalence}
As discussed above, it was shown in \cite{IR1} that for $3$-CY-algebras quiver mutation can be categorified using mutation of tilting modules of projective dimension at most $1,$ similar to the original case of categorifying reflections at sinks discussed in Section \ref{sec1}. And a tilting module gives rise to a derived equivalence \cite{Ha2}.

It is known from \cite{DWZ1} that for any finite quiver $Q$ without loops or $2$-cycles, there is some potential $W$ with the following property. For any quiver with potential $(Q',W')$ obtained from $(Q,W)$ by a finite sequence of mutations, the quiver $Q'$ has no loops or $2$-cycles. Then the operation $\mu_i$ taking $(Q,W)$ to $\mu_i(Q,W) = (\mu_i(Q),\overline{W})$ is directly a categorification of the quiver mutation taking $Q$ to $\mu_i(Q),$  without having to remove any $2$-cycles after performing the mutation $\mu_i$ on $(Q,W).$  The following result from \cite{KY} is a generalization of the results in  \cite{IR1}.

\begin{theorem}
With the above notation, the Ginzburg algebras associated with $(Q,W)$ and $\mu_i(Q,W)$ are derived equivalent.
\end{theorem}

\section{Applications to Cluster Algebras}\label{sec6}
The categorification of various classes of cluster algebras is an interesting problem itself. We have seen that the special case of categorifying quiver mutation in the acyclic case led to information on cluster algebras, namely a characterization of the acyclic cluster algebras having only a finite number of quivers occurring in the seeds. Categorification has also been used to discover new cluster algebras and to categorify old ones. In order to use categorification to obtain additional information on cluster algebras, it is of interest to define maps with nice properties between cluster variables and indecomposable rigid objects, and show that they are injective and/or surjective.

\subsection{The Dynkin case}
For Dynkin diagrams it was already known from \cite{FZ3} and \cite{MRZ} that there is a bijection between the cluster variables and the almost positive roots, hence a bijection between cluster variables and indecomposable decorated representations. For the associated cluster category $\CQ$ we have a natural correspondence between the negative simple roots $\seqs{-s}$ and the indecomposable objects $P_1[1],\ldots,P_n[1]$ in the cluster category.  Here $P_i$ is the projective cover of the simple $kQ$-module $S_i$ corresponding to the simple root $s_i.$ Then we have the following \cite{BMRRT}.

\begin{theorem}\label{Thm61}
Let $Q$ be a Dynkin quiver of type $A_n,\,D_n,\,E_6,\,E_7$ or $E_8$ and $\mathcal{A}(Q)$ the associated cluster algebra with no coefficients. Then there is a bijection from the cluster variables of $\mathcal{A}(Q)$ to the indecomposable (rigid) objects in the cluster category $\CQ,$ sending clusters to cluster tilting objects.
\end{theorem}

\subsection{From cluster algebras to cluster categories}\label{sec62}
It was conjectured in \cite{BMRRT} that there should be a bijection as in Theorem \ref{Thm61} in the general acyclic case. Relevant maps have been defined in both directions in order to deal with this problem. Here we start with defining a natural map $\a$ from cluster variables to indecomposable rigid objects \cite{BMRT}.

Let $(\seqn{x},Q),$ with $Q$ a finite connected acyclic quiver, be the initial seed. Then define $\a(x_i) = P_i[1]$ for $i = \num{n}.$ The mutation $\mu_i(\seqn{x},Q)$ creates a new cluster variable $ x_i^*,$ which is sent to the indecomposable rigid object $P_i[i]^*,$ where $P_i[1]^*\nsimeq P_i[1]$ and $\left(kQ/P_i[1]\right)\oplus{P_i[1]^*}$ is a cluster tilting object in $\CQ.$ We continue this procedure, and prove that the map $\a$ is well-defined. Here we use heavily property (C4) of a cluster structure, as proved in \cite{BMR2}. To show that $\a$ is surjective one uses that the cluster tilting  graph is connected. Then one shows that the map $\a$ on cluster variables induces a map from clusters to cluster tilting objects.

Already from these properties of $\a$ one gets the following, which was proved independently in \cite{BMRT} and  \cite{CK2} for acyclic cluster algebras.

\begin{theorem}
For acyclic cluster algebras with no coefficients, a seed is determined by its cluster.
\end{theorem}

This has later been generalized in \cite{GSV} using other methods.

\subsection{Interpretation of denominators}
Let $(\seqn{x},Q)$ be the initial seed, where $Q$ is an acyclic quiver. The denominator of a cluster variable (different from $\seqn{x}),$ expressed in the variables $\seqs{x},$  determines the composition factors of a unique indecomposable rigid $kQ$-module. This is proved in \cite{BMRT} at the same time as constructing the surjective map $\a$ discussed in \ref{sec62}. For example, for $Q:$ 
\begin{minipage}{2.4cm}
\begin{center} 
\begin{tikzpicture}
\node[] () at (0,0) {\s $1$};
\node[] () at (1,0) {\s $2$};
\node[] () at (2,0) {\s $3$};
\draw[->] (0.1,0) -- (0.9,0);
\draw[->] (1.1,0) -- (1.9,0);
\end{tikzpicture}
\end{center}
\end{minipage} the cluster variable $f/x_1x_2x_3$ in reduced form is sent by $\a$ to the unique indecomposable rigid $kQ$-module which has composition factors $S_1,\,S_2,\,S_3,$ and this is $P_1.$

The surprisingly simple, but extremely useful, idea of
\textit{positivity condition} was crucial for the proof. This says that if $f = f\seqn{x}$ has the property that if $f(e_i)>0,$ where $e_i = (1,\ldots,1,\underset{i}{0},1,\ldots,1),$ for $i= \num{n},$ then $f/m,$ where $m$ is a non constant monomial, is in  reduced form.

The approach sketched above is taken from \cite{BMRT}. This type of connection between denominators and indecomposable rigid objects was first shown in \cite{FZ1} for Dynkin diagrams with alternating orientation, then in \cite{CCS2} and \cite{RT} for the general Dynkin case, with another approach in \cite{CK1}. Note that only using the approach sketched above there could still theoretically be  different cluster variables $f/m$ and $f'/m$ in reduced form, with the same monomial $m.$ Another approach to the denominator theorem is given in \cite{CK2}, where also the above positivity condition from \cite{BMRT} is used, together with the Caldero-Chapoton map, which we discuss next. Using this map, it follows that $m$ determines $f/m$ (see also \cite{Hub}).

Note that when we express the cluster variables in terms of a seed different from the initial one, the denominators do not necessarily determine an indecomposable rigid module \cite{BMR4}. This fact was useful in \cite{FK} for giving a counterexample to a conjecture in \cite{FZ5}.

\subsection{From cluster categories to cluster algebras}
We now define a map from indecomposable rigid objects in cluster categories to cluster variables in the corresponding cluster algebras. This beautiful work was started in \cite{CC} for the Dynkin case, with the following definition.

Let $M$ be in $\m kQ.$ Then define
$$X_M =\sum\limits_{\e}\chi(Gr_{\e}(M))\underset{i}{\prod}u_i^{-\langle\e,\a_i\rangle-\langle\a_i,\mathbf{m}-\e\rangle}.$$
Here $Gr_{\e}(M) = \{N\in\m kQ; N\subset M,\, \underline{\mbox{dim}}N
= \e\},$ where $\e = \seqn{e}\leq \mathbf{m} = \seqn{m},$ which denotes the dimension vector of $M.$ Further $\langle\,,\,\rangle$ denotes the Euler form defined on the  Grothendieck group of $\m kQ.$ For $M, N$ in $\m kQ $ it is defined by $\langle M,N\rangle = \mbox{dim}_k\Hom_{kQ}(M,N)-\mbox{dim}_k\ext_{kQ}^1(M,N).$ Finally $\chi$ denotes the Euler-Poincar\'e characteristic of the quiver Grassmanian. Further we define $X_{P_i[1]} = u_i$, where $P_1, P_2,\ldots, P_n$ are the
indecomposable projective $kQ$-modules. 

We illustrate with the following.

\begin{example}\normalfont
Let $Q$ be the quiver 
\begin{minipage}{2.3cm}
\begin{center} 
\begin{tikzpicture}
\node[] () at (0,0) {\s $1$};
\node[] () at (1,0) {\s $2$};
\node[] () at (2,0) {\s $3$};
\draw[->] (0.1,0) -- (0.9,0);
\draw[->] (1.1,0) -- (1.9,0);
\end{tikzpicture}
\end{center}
\end{minipage}. We compute $X_{S_2}.$ We have $\e = (0,0,0)$ or $\e = (0,1,0).$ In both cases $\chi(Gr_{\e}(S_2))=1.$

Assume first $\e = (0,0,0).$ Then we have 
$$u_1^{-\langle 0,\a_1\rangle -\langle \a_1,\a_2\rangle }u_2^{-\langle 0,\a_2\rangle -\langle \a_2,\a_2\rangle }u_3^{-\langle 0,\a_2\rangle -\langle \a_3,\a_2\rangle } = u_1/u_2$$
since $\ext_{kQ}^1(S_1,S_2)\simeq k$ and $\ext_{kQ}^1(S_3,S_2) = 0.$

Assume then that $\e = \a_2.$ Then we have 
$$u_1^{-\langle \a_2,\a_1\rangle }u_2^{-\langle \a_2,\a_2\rangle }u_3^{-\langle \a_2,\a_3\rangle } = u_3/u_2$$ Hence we get
$$X_{S_2} = u/u_2 + u_3/u_2 = (u_1+u_3)/u_2.$$
\end{example}

The following is the main result of \cite{CC}.

\begin{theorem}
In the Dynkin case, the $\mathbb{Q}[\seqs{u}]$-submodule of $\mathbb{Q}\seqn{u}$ generated by the $X_M$ for $M$ an indecomposable $kQ$-module, coincides with the cluster algebra $\mathcal{A}(Q),$ and the set of cluster variables is given by  $$\{u_i; 1\leq i\leq n\}\cup\{X_M; M \mbox{ indecomposable module in } \m{kQ}\}.$$
\end{theorem}

The Caldero-Chapoton formula was generalized to the acyclic case in \cite{CK2}. The corresponding map $\b$ from the indecomposable rigid objects in the cluster category to the cluster variables in the associated cluster algebra was shown to be surjective, and also injective by using the positivity condition from the previous subsection. One then has the following \cite{CK2} (see also \cite[Appendix]{BMRT}).

\begin{theorem}
There is a map $\b$ from the indecomposable objects in the cluster category $\CQ$ for a finite acyclic quiver $Q$ to the associated cluster algebra, $\mathcal{A}(Q),$ with no coefficients, such that there is induced
\begin{itemize}
	\item[(1)] a bijection from the indecomposable rigid objects to the cluster variables
	\item[(2)] a bijection between the cluster tilting objects in $\CQ$ and the clusters for $\mathcal{A}(Q).$
\end{itemize}
\end{theorem}

Some further feedback to acyclic cluster algebras is then given as a result of this tighter connection, for example the following \cite[Appendix]{BMRT}.

\begin{theorem}
If $\seq{u}$ is a cluster in an acyclic cluster algebra with no coefficients, there is for each $i = \num{n},$ a unique element $u_i^*\neq u_i$ such that $\{u_1,\ldots,u_i^*,u_{i+1},\ldots,u_n\}$ is a cluster.
\end{theorem}

This has later been generalized to cluster algebras of \textit{geometric} type in \cite{GSV}, using different methods.

There is another variation of the Caldero-Chapoton formula in
\cite{P}, where the case of an arbitrary initial seed in an acyclic
cluster algebra was treated. In that connection Palu formulated the desired properties needed in order to obtain a good map from a 2-CY triangulated category $\C$ to a commutative ring $R$, called a \textit{cluster character} $\chi$ \cite{P} (see also \cite{BIRSc}). The requirement was that 
\begin{itemize}
	\item[(i)] $\chi(A) = \chi(B)$ if $A\simeq B$.
	\item[(ii)] $\chi(A\oplus B) = \chi(A)\chi(B)$.
	\item[(iii)] If $\mbox{dim}\ext^1(X,Y)=1$ for indecomposable
	objects $X$ and $Y$ in $\C$, consider the non split triangles $X\rightarrow B\rightarrow Y\rightarrow X[1]$ and $Y\rightarrow B'\rightarrow X\rightarrow Y[1].$ Then we have  $\chi(X)\chi(Y) = \chi(B) + \chi(B').$
\end{itemize}

The last condition is needed to ensure that indecomposable rigid objects are sent to cluster variables (see also \cite{BIRSc}).

There are several interesting generalizations beyond the acyclic case, and we refer to (\cite{DK},\cite{FK},\cite{GLS6},\cite{GLS1},\cite{GLS5},\cite{Nag},\cite{Pl1},\cite{Pl2}).

\section*{Acknowledgement}
I want to express my gratitude to Aslak Bakke Buan, Osamu Iyama, Bernhard Keller and Robert Marsh for helpful comments and suggestions, and I would like to thank them and my other coauthors on the topic of this paper for enjoyable collaboration.

\end{document}